\documentclass[11pt]{amsart}
\usepackage[utf8]{inputenc}
\usepackage[english]{babel}
\usepackage{romanbar}
\usepackage{amscd,amssymb}
\usepackage{amsmath}
\usepackage{amsthm}
\usepackage{graphicx}
\usepackage{fancybox}
\usepackage{epic,eepic}
\usepackage{amstext}
\usepackage{dsfont}
\usepackage{esint}
\usepackage[square,numbers]{natbib}
\usepackage{booktabs}
\usepackage{hyperref}
\usepackage[noabbrev, capitalise, nameinlink]{cleveref}
\crefname{equation}{}{}
\crefname{assumption}{Assumption}{Assumptions}
\crefformat{equation}{\textup{#2(#1)#3}}
\usepackage{comment}
\usepackage{mathtools}
\usepackage{tikz,pgfplots}
\usepackage{subcaption}
\usepackage{easytable}

\def\Nb{{\mathsf N}}
\newcommand{\with}{\,\colon\,}

\def\R{\mathcal{R}}

\def\Pnull{\mathbb{P}^0}

\def\R{{\mathbb R}}

\def\ba{\boldsymbol{a}}
\def\bb{\boldsymbol{b}}

\def\bA{\boldsymbol{A}}

\def\bv{\boldsymbol{v}}

\def\bw{\boldsymbol{w}}
\def\bg{\boldsymbol{g}}

\def\bvarphi{\boldsymbol{\varphi}}

\def\AA{{\mathcal A}}
\def\BB{{\mathcal B}}

\def\TT{{\mathcal T}}

\def\AAb{\boldsymbol{\AA}}

\def\ub{\boldsymbol{u}}
\def\vb{\boldsymbol{v}}

\def\Ab{\boldsymbol{A}}
\def\ab{\boldsymbol{a}}
\def\bb{\boldsymbol{b}}

\def\Yb{\boldsymbol{Y}}

\def\E{{\mathbb E}}

\def\norm#1#2{\|#1\|_{#2}}

\def\id{\mathds{1}}

\def\normL2#1#2{\|#1\|_{L^2(#2)}}
\def\normH10#1#2{\|#1\|_{H^1_0(#2)}}

\def\T{\mathcal{T}}
\def\TH{\mathcal{T}_H}

\def\PiH {\Pi_H}

\def\usbar{\bar{u}_{H,\ell}}
\def\uH{\ub_{H,\ell}}

\def\philoc{\phideal^\mathrm{loc}}
\def\phideal{\bvarphi_T}

\def\norml2l2#1#2{\|#1\|_{L^2(\Omega,L^2(#2))}}

\def\corrc{\boldsymbol{\mathcal{C}}_{T,\ell}}
\def\philod{\bvarphi_{T,\ell}^\mathrm{LOD}}
\def\glod{\bg_{T,\ell}^\mathrm{LOD}}
\def\ptl{\boldsymbol{p}_{T, \ell}}
\def\qtl{\boldsymbol{q}_{T, \ell}}

\def\Apatch{\AAb_{T, \ell}}
\def\Apatchinv{\AAb^{-1}_{T, \ell}}

\def\tr{{\operatorname{tr}}}
\def\trinv{\boldsymbol{\operatorname{tr}}^{-1}}

\def\Crb{C_\mathrm{rb}}

\newcommand{\tnormf}[2]{{\| #1 \|}_{L^2(#2)}}

\newcommand{\tsp}[3]{{\left( #1\,,\,#2 \right)}_{L^2(#3)}}
\newcommand{\tspf}[3]{{( #1\,,\,#2 )}_{L^2(#3)}}

\newcommand\dx{\,\mathrm{d}x}

\newcommand\dtx{\,\mathrm{d}\tilde x}

\definecolor{myBlue}{RGB}{113,104,238} 
\definecolor{myGreen}{RGB}{114,175,30} 
\definecolor{myRed}{RGB}{180,50,50}  
\definecolor{myOrange}{RGB}{225,92,22} 

\newtheorem{theorem}{Theorem}[section]
\newtheorem{corollary}[theorem]{Corollary}
\newtheorem{lemma}[theorem]{Lemma}
\newtheorem{assumption}[theorem]{Assumption}

\theoremstyle{definition}

\theoremstyle{remark}
\newtheorem{remark}[theorem]{Remark}
\numberwithin{theorem}{section}
\numberwithin{equation}{section}
\numberwithin{table}{section}
\numberwithin{figure}{section}

\allowdisplaybreaks

\begin{document}
	
\title[Numerical stochastic homogenization]{A simple collocation-type approach to numerical stochastic homogenization}
\author[M.~Hauck, H.~Mohr, D.~Peterseim]{Moritz Hauck$^\dagger$, Hannah Mohr$^\ddagger$, Daniel Peterseim$^{\ddagger,\star}$}
\address{${}^{\dagger}$ Interdisciplinary Center for Scientific Computing (IWR), Heidelberg University, 69047 Heidelberg, Germany}
\email{moritz.hauck@iwr.uni-heidelberg.de}
\address{${}^{\ddagger}$ Institute of Mathematics, University of Augsburg, Universit\"atsstr.~12a, 86159 Augsburg, Germany}
\address{${}^{\star}$ Centre for Advanced Analytics and Predictive Sciences (CAAPS), University of Augsburg,
	Universit\"atsstr.~12a, 86159 Augsburg, Germany}
\email{\{hannah.mohr, daniel.peterseim\}@uni-a.de}
\thanks{The work of Moritz Hauck, Hannah Mohr, and Daniel Peterseim is part of a project that has received funding from the European Research Council (ERC) under the European Union's Horizon 2020 research and innovation programme (Grant agreement No.~865751 --  RandomMultiScales). The work of Moritz Hauck is also supported by the Knut and Alice Wallenberg foundation postdoctoral program in mathematics for researchers from outside Sweden (Grant No.~KAW 2022.0260).}
\maketitle

\begin{abstract}
This paper proposes a novel collocation-type numerical stochastic homogenization method for prototypical stochastic homogenization problems with random coefficient fields of small correlation lengths. The presented method is based on a recently introduced localization technique that enforces a super-exponential decay of the basis functions relative to the underlying coarse mesh, resulting in considerable computational savings during the sampling phase. More generally, the collocation-type structure offers a particularly simple and computationally efficient construction in the stochastic setting with minimized communication between the patches where the basis functions of the method are computed. An error analysis that bridges numerical homogenization and the quantitative theory of stochastic homogenization is performed. In a series of numerical experiments, we  study the effect of the correlation length and the discretization parameters on the approximation quality of the method. 
\end{abstract}
\vspace{1cm}
\noindent\textbf{Keywords:}
numerical homogenization,
stochastic homogenization,
super-localization,
quantitative theory,
error estimates,
uncertainty
\\[2ex]
\textbf{AMS subject classifications:}
35R60,		
65N12, 		
65N15, 		
65N35,		
73B27		

\section{Introduction}
This paper presents a novel numerical stochastic homogenization method for the prototypical random diffusion problem
\begin{equation*}
	-\mathrm{div}(\Ab \nabla \ub) = f
\end{equation*} 
subject to homogeneous Dirichlet boundary conditions. Microscopic features of the problem are encapsulated in the random diffusion coefficient $\Ab$. In this paper, we are particularly interested in coefficients with small correlation lengths.

%

For deterministic coefficients, numerical homogenization techniques have been studied extensively in the last decades. For recent monographs and reviews on this topic, we refer to \cite{CEH23, BLeB23, AHP21, MalP20, OS19}. The random case has not received a similar attention. However, there are several numerical approaches. Let us mention, for example, MsFEM-based numerical stochastic homogenization methods that assume so-called weakly random coefficients; see~\cite{ACL11} for an overview. A popular approach to approximate the effective coefficient in stochastic homogenization is the so-called Representative Volume Element (RVE) method. Its theoretical analysis was first achieved in \cite{Glo12} in the discrete random setting and generalized in \cite{Gloria2014, GloNo16, GloHa16}; see also the works \cite{BoP04,Cancs2015,Mourrat2018,Fischer2019,Khoromskaia2020}.
Closely related to the present work is the numerical stochastic homogenization method proposed and analyzed in \cite{GaP19,FisGP19ppt}. This method is based on the Localized Orthogonal Decomposition (LOD) introduced in~\cite{MaP14,HeP13}; see also~\cite{MalP20,AHP21} for an overview. More precisely, using a reformulation of the LOD for deterministic problems based on a quasi-local discrete integral operator, as discussed in \cite{GaP17}, one can derive an effective model of the problem at hand by taking the expectation. This effective model is deterministic, and its solution gives an accurate coarse scale approximation to the expected value of the solution. However, the method is strongly tied to linear finite elements on simplicial meshes (with piecewise constant gradients), which seems to be an artificial limitation and practically unfavorable, especially since structured deterministic and random diffusion coefficients are often based on Cartesian meshes. 

In addition to this technical shortcoming, the localization technique underlying the LOD framework has recently been improved, evolving into the Super-Localized Orthogonal Decomposition (SLOD) introduced in~\cite{HaPe21b} (see also \cite{Freese-Hauck-Peterseim,BHP22,Bonizzoni-Freese-Peterseim,pumslod,graphSLOD}).  While the LOD has exponentially decaying basis functions that lead to an exponentially decaying localization error with respect to the diameter of the basis supports relative to the underlying mesh, the localization error of the SLOD actually decays super-exponentially, cf.~\cite{HaPe21b}. This results in smaller local fine-scale problems when computing the basis functions and increased sparsity for a given tolerance. 

In this paper, we propose a computationally simple and  efficient numerical stochastic homogenization method based on a special collocation-type formulation of the SLOD. This collocation-type formulation leads to a coarse stiffness-type matrix that can be assembled without any communication between the basis functions defined on the patches of the coarse mesh. This allows each patch to be considered separately for sampling, allowing for improved parallelization and a significant speed-up of the method's assembly process. Furthermore, the favorable localization properties of the SLOD allow for a computationally efficient sampling procedure.
In the case of a random diffusion coefficient with a small correlation length and under standard assumptions of quantitative stochastic homogenization, this paper provides an error estimate for the coarse scale approximation of the proposed method, where certain SLOD-specific quantities contribute in an a posteriori manner. The proof of this error estimate is based on the theory of quantitative stochastic homogenization; see, e.g.,~\cite{Gloria2011,Gloria2012,Gloria2014,GloriaNeukammOttoPreprint}. Classical LOD-techniques \cite{MaP14,HeP13,AHP21} are used to further evaluate these SLOD-specific quantities, and a worst-case a priori error analysis is conducted for one of them. Several numerical experiments are performed to quantitatively study the effect of the correlation length and other discretization parameters on the accuracy of the approximation. 

This manuscript is structured as follows.  First, in \cref{sec:modelproblem} we state the model problem in weak form. Then, in \cref{sec:method}, we introduce the novel numerical stochastic homogenization method. An a posteriori error analysis of the method is performed in \cref{sec:erroranalysis}. \cref{sec:erroranalysiswithlod} derives a worst-case estimate for the quantity appearing in the a posteriori error bound. Practical aspects of implementation are addressed in \cref{sec:implementation}. Finally, \cref{sec:numexp} provides numerical experiments that underline the theoretical results of this paper.

\section{Model problem}\label{sec:modelproblem}
We consider the model problem
\begin{equation}\label{e:modelstrong}
	\left\{
	\begin{aligned}
			-\mathrm{div} (\Ab(\omega)(x) \nabla \ub(\omega)(x) ) &= f(x),&x\in D\\
			\ub(\omega)(x)&=0,&x\in \partial D
		\end{aligned}
	\right\}\quad\text{for almost all }\omega\in\Omega,
\end{equation}
where $(\Omega,\mathcal{F},\mathbb{P})$ denotes the underlying probability space, $f \in L^2(D)$ is a deterministic right-hand side and $D$ is a $d$-dimensional bounded Lipschitz polytope with $d\in\{1,2,3\}$. Without loss of generality, we assume that $D$ is scaled to unit size. 
Suppose that $\Ab$ is a $\mathbb R^{d\times d }$-valued pointwise symmetric Bochner measurable function, which is uniformly elliptic and bounded, i.e., there exist $0<\alpha \leq \beta < \infty$ such that for almost all $\omega \in \Omega$ 
\begin{equation}
	\label{eq:coeffspd}
	\alpha |\xi|^2 \leq \langle \xi,\Ab(\omega)(x)\xi \rangle \leq \beta |\xi|^2	
\end{equation}
holds for all $\xi \in \mathbb R^d$ and almost all $x \in \mathbb R^d$, where $\langle\cdot,\cdot\rangle$ denotes the Euclidean inner product of $\mathbb R^d$ and $|\cdot|$ its induced norm. Note that the above symmetry assumption on $\Ab$ is made for the sake of simplicity. In fact, we expect that the construction of the proposed method and the corresponding proofs can be easily generalized to the non-symmetric case.

The weak formulation of the model problem~\cref{e:modelstrong} seeks a $H^1_0(D)$-valued random field~$\ub$ such that for almost all $\omega\in\Omega$ it holds that  
\begin{equation}\label{e:modelweak}
\ab_\omega(\ub(\omega),v):=\tspf{\Ab(\omega)\nabla\ub(\omega)}{\nabla v}{D}= \tspf{f}{v}{D}\quad\text{for all }v\in H^1_0(D).
\end{equation}
Here, $\tsp{\cdot}{\cdot}{D}$ denotes the inner product on $L^2(D)$ or $(L^2(D))^d$.

Subsequently, we introduce a shorthand notation for norms and inner products of Bochner spaces. Let $X$ be a Hilbert space equipped with the inner product $(\cdot,\cdot)_X$. In this case, the Bochner space $L^2(\Omega;X)$, denoting the space of $X$-valued random fields with finite second moments, is also a Hilbert space with the inner product 
\begin{equation*}
	(\bv,\bw)_{L^2(\Omega;X)}\coloneqq \E\big[(\bv(\omega),\bw(\omega))_{X}\big].
\end{equation*}
We write $\|\cdot\|_{L^2(\Omega;X)}^2 \coloneqq (\cdot,\cdot)_{L^2(\Omega;X)}$ for the induced norm of this inner product. 

Under the given assumptions, the bilinear from~$\ab_\omega$ depends continuously on $\Ab$ and, in particular, is measurable as a function of~$\omega$. Hence, the above problem can be reformulated in the Hilbert space $L^2(\Omega;H^1_0(D))$. The Lax--Milgram theorem then proves its well-posedness, i.e., there exists a unique solution ${\ub\in L^2(\Omega;H^1_0(D))}$ satisfying
\begin{equation}
	\label{eq:stabilitysol}
	\|\nabla\ub\|_{L^2(\Omega;L^2(D))} \leq \alpha^{-1}C_{\mathrm F}\|f\|_{L^2(D)},
\end{equation}
where $\tnormf{\cdot}{D}$ denotes the $L^2(D)$-norm and $C_F$ is the Friedrichs constant of~$D$.

\section{Numerical stochastic homogenization method}\label{sec:method}

The construction of the novel stochastic homogenization method is based on ideas of the SLOD introduced in \cite{HaPe21b}. In the deterministic setting, the SLOD identifies an almost local basis of the space obtained by applying the solution operator to $\Pnull(\TH)$, the space of piecewise constants with respect to some coarse mesh $\TH$ of the domain $D$. This is achieved by identifying piecewise constant right-hand sides supported on patches of the mesh~$\TH$ such that their responses under the corresponding localized solution operator have a minimal conormal derivative. These localized responses are then used as the basis functions of the SLOD. In the stochastic setting, an adaptation of this approach is required, which involves identifying deterministic local source terms such that the conormal derivative of the localized responses is small in expectation. 

In the following we assume that the considered family of meshes $\{\TH\}_{H}$ is quasi-uniform and consists of meshes with closed, convex and shape-regular elements. The parameter $H>0$ specifies the maximal element diameter of the mesh $\mathcal T_H$. We denote by $\Pi_H\colon L^2(D) \rightarrow \mathbb{P}^0(\TT_H)$ the $L^2$-orthogonal projection onto $\mathbb{P}^0(\TT_H)$. Let us also give a precise definition of the concept of patches with respect to~$\TH$. The first-order patch $\Nb(S) = \Nb^1(S)$ of $S \subset D$, where $S$ is a union of elements of $\TH$, is defined by 
\begin{equation*}	
	\Nb^1(S) \coloneqq \bigcup \{T\in\TT_H \with T\cap S \neq \emptyset \}.
\end{equation*}
For any $\ell = 2,3,4,\dots$, the $\ell$-th order patch $\Nb^\ell(T)$ of $T$ is then given recursively by 
\begin{equation}
	\label{eq:patch}
	\Nb^\ell(T) \coloneqq \Nb^1(\Nb^{\ell-1}(T)).
\end{equation}

The following derivation of the basis functions considers a fixed element $T\in \TT_H$ and oversampling parameter $\ell \in \mathbb N$, where we assume that the  patch $D_T \coloneqq \Nb^\ell(T)$ does not coincide with the whole domain. We denote the deterministic source term corresponding to $T$ by $g_T \in \Pnull (\TT_{H,D_T})$ with the submesh $\TT_{H,D_T} \coloneqq \{K \in \TT_H \with K \subset D_T\}$. Note that in the following we do not distinguish between locally defined $L^2$- or $H^1_0$-functions and their  extensions by zero to the whole domain. The global response $\phideal \in L^2(\Omega; H^1_0(D))$ to $g_T$ is then defined for almost all $\omega \in \Omega$ by
\begin{equation}
	\label{eq:basisglob}
	\ab_\omega(\phideal(\omega),v) = \tspf{g_T}{v}{D}\quad\text{for all }v\in H^1_0(D).
\end{equation}
Its localized version $\philoc \in L^2(\Omega; H^1_0(D_T))$ is for almost all $\omega \in \Omega$ defined by
\begin{equation}
	\label{eq:basisloc}
	\ab_\omega(\philoc(\omega),v) = \tspf{g_T}{v}{D_T}\quad\text{for all }v\in H^1_0(D_T).
\end{equation}
From now on, the dependence of stochastic variables is only indicated by a bold symbol for better readability.

 To define the conormal derivative of the localized basis function $\philoc$, we need to introduce some preliminaries on traces and extensions. We denote by $H^1_\Gamma(D_T)$ the complete subspace of $H^1(D_T)$ consisting of functions with trace zero at the boundary segment $\Gamma \coloneqq \partial D_T \cap \partial D$. Furthermore, let
\begin{align*}
	\tr\colon H^1_\Gamma(D_T)\to X \coloneqq \operatorname {im} \tr  \subset H^{1/2}(\partial D_T)
\end{align*}
denote the classical trace operator restricted to $H^1_\Gamma(D_T)$. As an extension operator, we henceforth consider the $\Ab$-harmonic extension  operator $\trinv\colon L^2(\Omega;X) \to L^2(\Omega; H^1_\Gamma(D_T))$ defined  as follows: For almost all $\omega \in \Omega$ and for any given $\bb\in L^2(\Omega;X)$, we set $(\tr\trinv \bb)(\omega) = \bb(\omega)$  and demand that
\begin{align}\label{eq:trinv}
	\ba (\trinv \bb, v ) = 0\quad \text{for all } v\in H^1_0(D_T).
\end{align}
The space of locally
$\Ab$-harmonic functions satisfying homogeneous Dirichlet boundary conditions on~$\Gamma$ can then be defined as
\begin{align*}
	\Yb \coloneqq \trinv L^2(\Omega;X) \subset L^2(\Omega; H^1_\Gamma(D_T)).
\end{align*}
For more details on trace- and extension operators we refer, e.g., to~\cite{LiM72a}.

Combining \cref{eq:basisloc,eq:trinv} yields for almost all $\omega \in \Omega$ and all $\bv \in L^2(\Omega;H^1_0(D))$ the identity
\begin{equation*}
	\ba(\philoc, \bv) = \ba(\philoc, \bv - \trinv \tr \bv) = (g_T, \bv - \trinv \tr \bv)_{L^2(D_T)},
\end{equation*}
where we used that $(\bv - \trinv \tr \bv)|_{D_T} \in L^2(\Omega;H^1_0(D_T))$. With this identity, the definition of $\phideal$ in \cref{eq:basisglob}, and $\operatorname{supp} g_T \subset D_T$, it follows that
\begin{align*}
	\ab(\phideal - \philoc, \bv) = (g_T, \bv)_{L^2(D_T)} - \ba(\philoc, \bv) = (g_T, \trinv \tr \bv)_{L^2(D_T)}.
\end{align*}
Taking the expectation, we obtain for any $\bv \in L^2(\Omega;H^1_0(D))$ that
\begin{align}\label{eq:loc-error}
	\E \big[\ab(\phideal - \philoc, \bv) \big] = \E\big[(g_T, \trinv \tr \bv)_{L^2(D_T)}\big] = \left(g_T, \Pi_{H, D_T}\; \E[ \trinv \tr \bv]\right)_{L^2(D_T)}.
\end{align}
As a consequence, the (almost) $L^2$-orthogonality of $g_T$ to the space $\E[\Yb] \subset H^1_\Gamma(D_T)$ leads to a small expected localization error for the basis function $\philoc$.

Therefore, we can obtain an optimal choice of $g_T$ by performing a singular value decomposition (SVD) of the compact operator $(\Pi_{H, D_T} \circ \E)\vert_{\Yb}: \Yb \rightarrow \mathbb P^0 (\TT_{H, D_T})$ restricted to the complete subspace $\bf Y$. Note that the rank of $(\Pi_{H, D_T} \circ \E)\vert_{\Yb}$ is less than or equal to $N \coloneqq \# \TT_{H,D_T}$. Hence, the SVD is given by 
\begin{equation}
	\label{eq:svd}
	(\Pi_{H, D_T} \circ \E)\vert_{\Yb} \; \bv = \sum_{k=1}^{N} \sigma_k (\bv,  \bw_k)_{L^2(\Omega; H^1(D_T))} g_k
\end{equation}
with singular values $\sigma_1 \geq \dots \geq\sigma_N\geq 0$, $L^2(\Omega; H^1(D_T))$-orthonormal right singular vectors $\bw_1,\dots,\bw_N$ and $L^2(D_T)$-orthonormal left singular vectors $g_1,\dots, g_N$. The choice $g_T=g_N$ as the left singular vector corresponding to the smallest singular value $\sigma_N$ is optimal in the sense that 
\begin{equation*}
	g_N \in \underset{g\in\mathbb{P}^0(\TT_{H,D_T})\;:\; \Vert g \Vert_{L^2(D_T)}=1}{\text{argmin}} \quad \underset{\bv\in \Yb\; :\; \Vert \bv \Vert_{L^2(\Omega;H^1(D_T))}=1}{\sup} (g,\E[\bv])_{L^2(D_T)}.
\end{equation*}
The corresponding smallest singular value $\sigma_N$ is a measure of the (quasi-)orthogonality between $g_T$ and $\E[\Yb]$. We hence define
\begin{equation}\label{eq: sigma}
	\sigma_T(H,\varepsilon,\ell) \coloneqq \sigma_N = \underset{\bv\in\Yb\; :\; \Vert \bv \Vert_{L^2(\Omega;H^1(D_T))}=1}{\sup} (g_N,\E[\bv])_{L^2(D_T)},
\end{equation}
where the parameter $\varepsilon>0$ denotes the correlation length of the random coefficient $\Ab$, which will be rigorously introduced in \cref{assumption_coefficient} below.

We emphasize that the practical implementation of the SVD in~\cref{eq:svd} is difficult due to the stochasticity involved; for a practical implementation based on sampling, see \cref{sec:implementation}. For the error analysis in the following section, we introduce the quantity
\begin{equation}
	\label{eq:locerrind}
	\sigma\coloneqq\sigma(H,\varepsilon,\ell) \coloneqq \max_{T\in \TH} \sigma_T(H,\varepsilon,\ell),
\end{equation}
which is an indicator for the  overall localization error.

\begin{figure}
	\includegraphics[width=.3\linewidth]{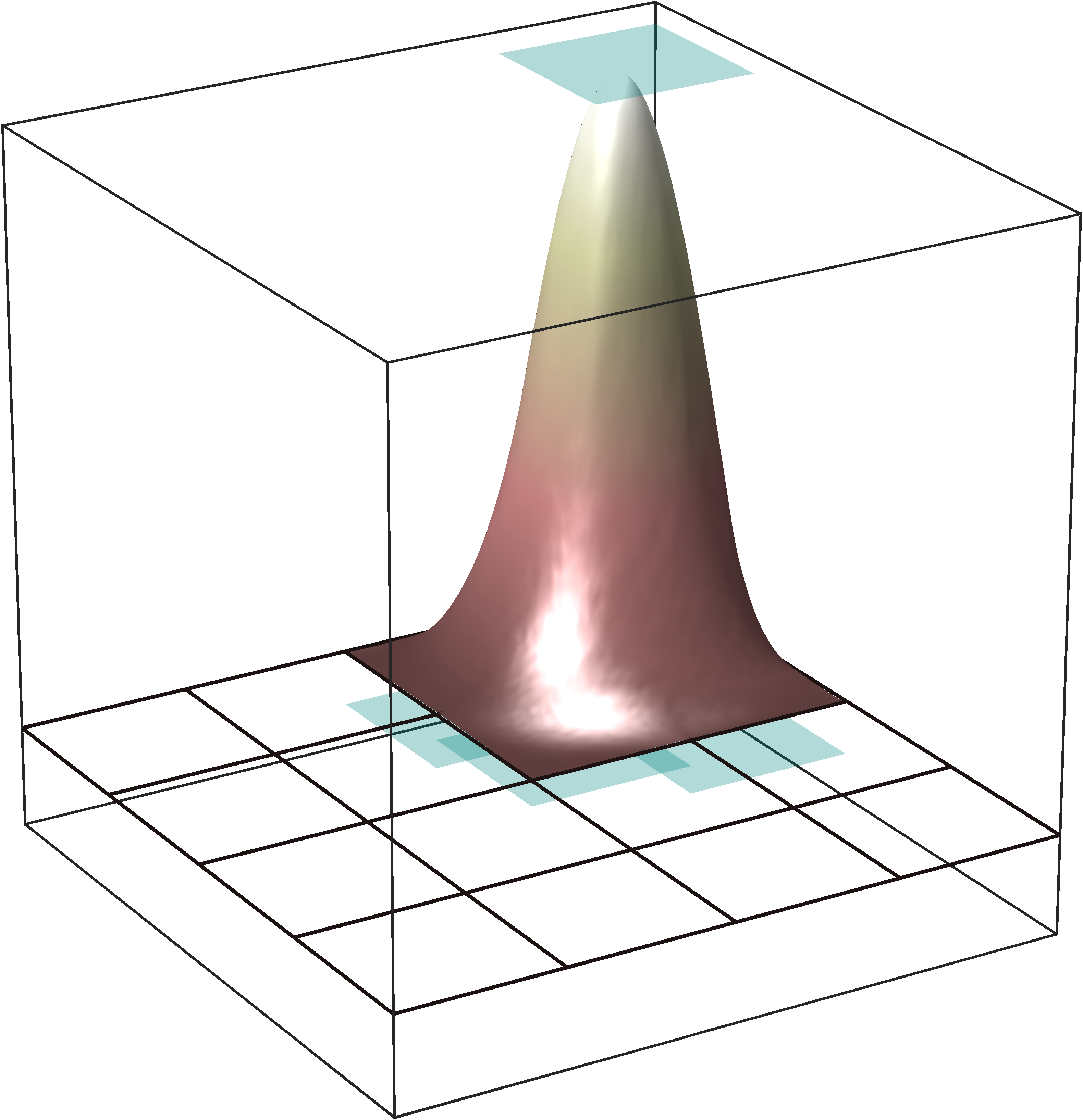}
	\hfill
	\includegraphics[width=.3\linewidth]{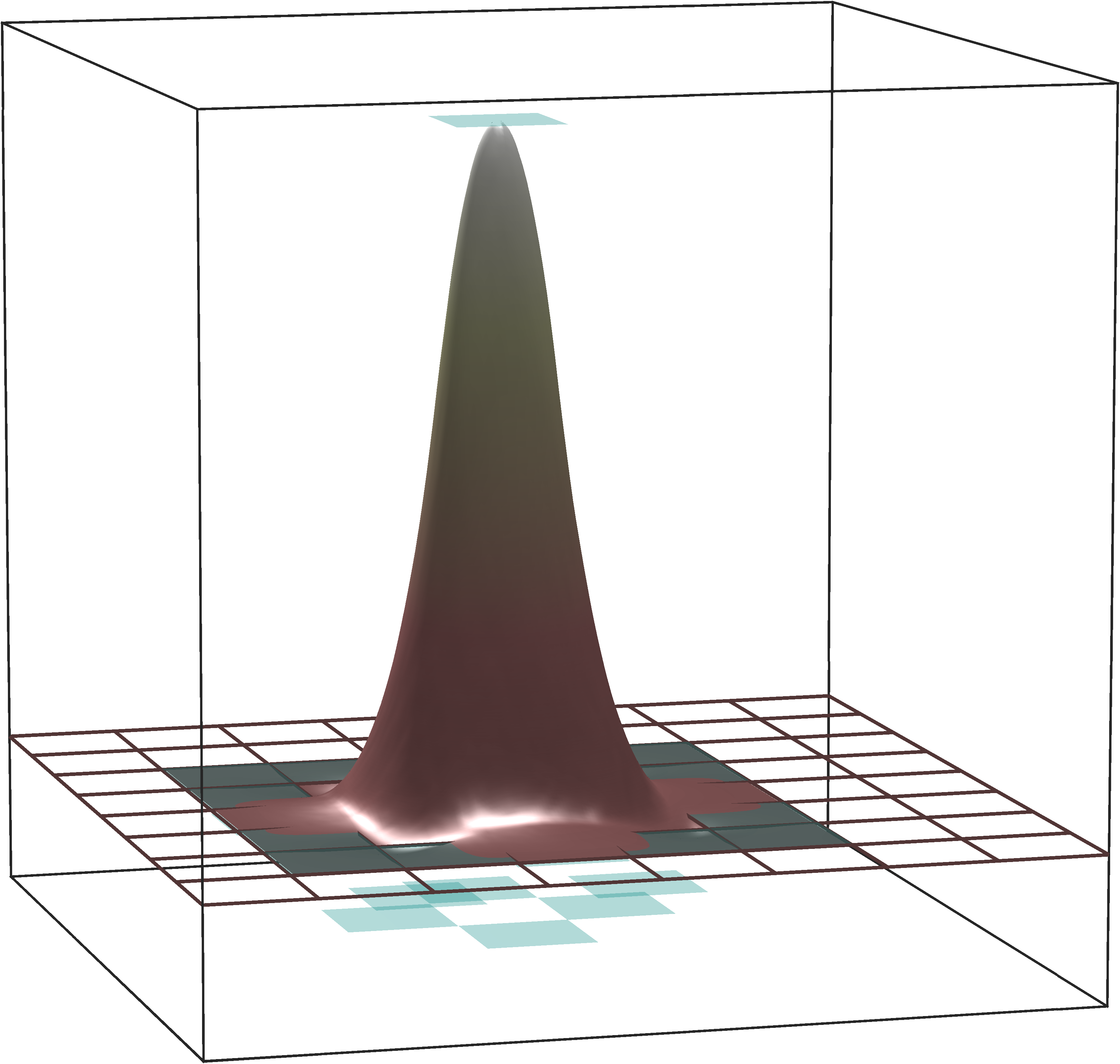}
	\hfill
	\includegraphics[width=.3\linewidth,height=0.275\linewidth]{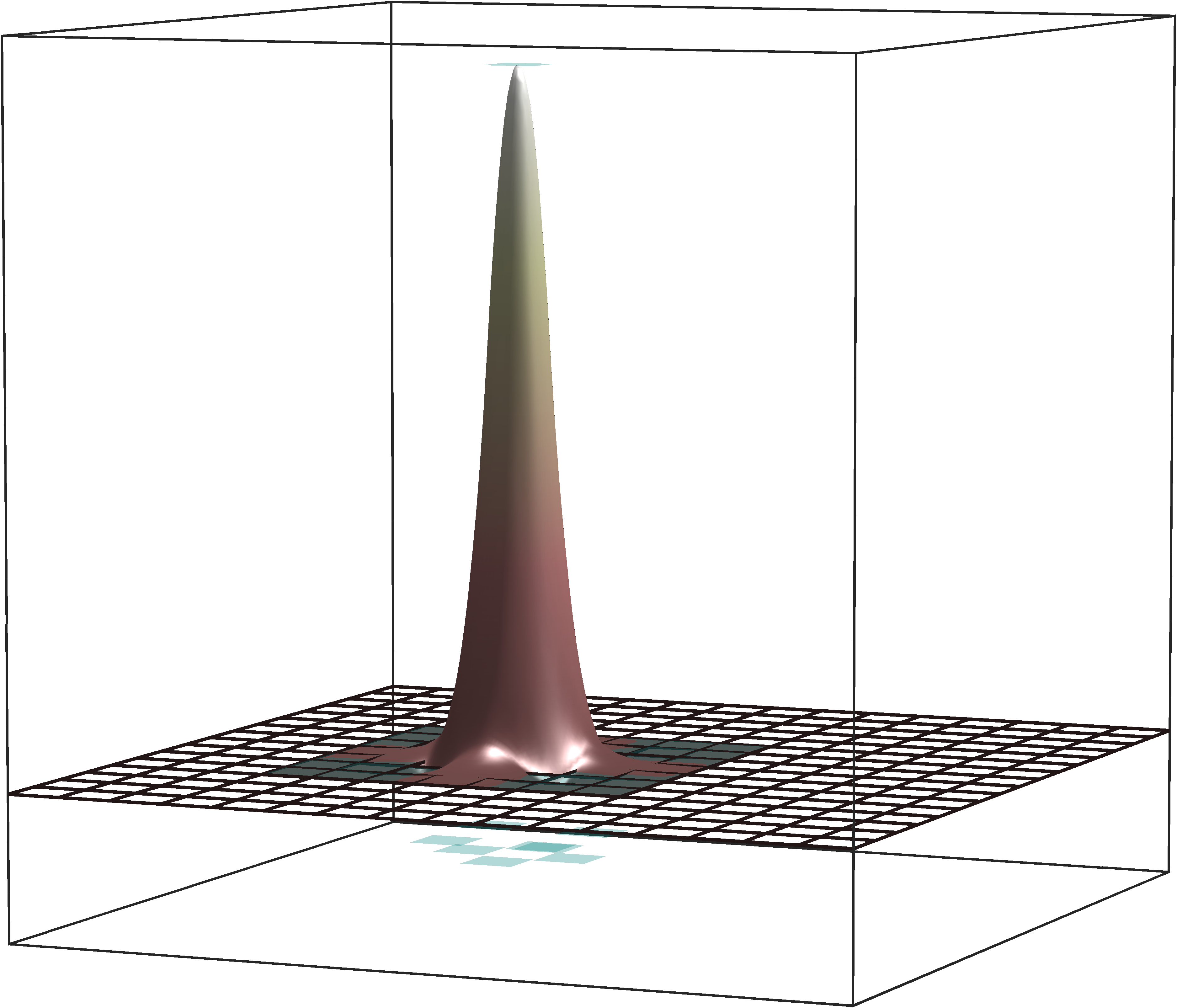}
	\caption{Illustration of the localized basis functions $\E[\boldsymbol{\varphi}^\mathrm{loc}]$ obtained by the novel stochastic homogenization method on successively refined meshes for a piecewise constant random coefficient with a correlation length of $\varepsilon = 2^{-7}$ in two spatial dimensions. Various values of the oversampling parameter are depicted with $\ell = 1$ (left), $\ell = 2$ (middle), and $\ell = 3$ (right). The corresponding right-hand sides $g$ are shown in green.}
	\label{fig:basis2d}
\end{figure}

Given that, in expectation, $\philoc$ closely approximates the response of the global solution operator applied to $g_T$, it is reasonable to define the approximation of a non-Galerkin, collocation-type numerical stochastic homogenization method by
\begin{equation}
	\label{eq:stochhommethod}
	\usbar \coloneqq \sum_{T\in\TT_{H}} c_T \PiH \mathbb{E}[\philoc],
\end{equation}
where $(c_T)_{T\in \TT_{H}}$ are the coefficients of the expansion of $\PiH f$ in terms of the basis functions $\{g_T \with T \in \TH\}$. An illustration of the deterministic basis functions $\E[\philoc]$ can be found in \cref{fig:basis2d}. Error estimates for this method are derived in the following sections. 

\section{Error analysis}
\label{sec:erroranalysis}

In this section, we perform an error analysis of the proposed stochastic homogenization method based on results from the theory of quantitative stochastic homogenization. This theory requires structural conditions on the randomness of the coefficient field $\Ab$. For simplicity, the conditions are formulated for coefficient fields defined on $\mathbb R^d$. Hence, the following assumptions implicitly assume that the coefficient field is defined on the full space~$\mathbb R^d$. A random field defined on the bounded domain $D$ can be obtained by restriction.

\begin{assumption}[Stationarity and decorrelation]\label{assumption_coefficient}
	Assume that the random coefficient field $\Ab$ is
	\begin{itemize}
		\item  \emph{stationary}, i.e., the law of the shifted coefficient field $\Ab(\omega)(\cdot+x)$ coincides with the law of $\Ab(\omega)(\cdot)$ for all $x\in \mathbb R^d$,
		\item \emph{quantitatively decorrelated} on scales larger than $\varepsilon$ in the sense of the spectral gap inequality with correlation length $\varepsilon>0$, i.e., there exists a constant $\rho>0$ such that for any Fr\'echet differentiable random variable $F=F(\Ab)$ the estimate
		\begin{align}
			\label{SpectralGap}
			\mathbb{E}\big[|F-\mathbb{E}[F]|^2\big]
			\leq \frac{\varepsilon^d}{\rho} \mathbb{E}\Bigg[\int_{\mathbb R^d} \bigg(\fint_{B_\varepsilon(x)} \bigg|\frac{\partial F}{\partial \Ab}(\tilde x)\bigg| \dtx\bigg)^2 \dx\Bigg]
		\end{align}
		holds.
	\end{itemize}
\end{assumption}
For an introduction to the notion of Fr\'echet derivatives, we refer the reader exemplarily to \cite[Chap.~2]{Deimling1985}; see also \cite[Sec. 3.1]{JO22} for a definition in the present context.  We emphasize that the conditions in \cref{assumption_coefficient} on the random coefficient $\Ab$ are standard in the theory of quantitative stochastic homogenization; see, e.g., the work \cite{GloriaNeukammOttoPreprint}.

The error bound  presented in this section is an a posteriori bound including the constant~$\sigma$ from \cref{eq:locerrind} and the Riesz stability constant of the local source terms $\left\{g_T\with T\in\TH\right\}$, which quantifies their linear independence. Both constants can be computed a posteriori as outlined in \cref{sec:implementation}. Additionally, we provide a worst-case a priori upper bound on $\sigma$ in \cref{sec:erroranalysiswithlod}. Note that in a practical implementation, the Riesz stability of the local source terms can be ensured as outlined in \cref{sec:implementation} or \cite[App.~B]{HaPe21b}.
\begin{assumption}[Riesz stability]\label{riesz_stability}
	The set $\left\{g_T\with T\in\TH\right\}$ is a Riesz basis of  $\mathbb{P}^0(\TH)$, i.e., there exists  $C_{\mathrm{rb}}(H, \ell)>0$ such that for all possible choices of $(c_T)_{T\in\TH}$ it holds
	\begin{equation*}
		C^{-1}_{\mathrm{rb}}(H,\ell) \sum_{T\in\TH} c^2_T \leq \bigg\| \sum_{T\in\TH} c_T g_T \bigg\|^2_{L^2(D)},
	\end{equation*}
where $C_\mathrm{rb}(H,\ell)$ 
depends polynomially on $H^{-1}$ and $\ell$.
\end{assumption}

To handle the stochasticity in the error analysis, we need to estimate the variance of the random variables $(\philoc, \id_K)_{L^2(K)}$ for any $T,K \in \TH$, where $\id_K$ denotes the indicator function of the element $K$. To achieve this, we employ the spectral gap inequality \cref{SpectralGap} from \cref{assumption_coefficient}. The following lemma provides a representation of the Fr\'echet derivative of $(\philoc, \id_K)_{L^2(K)}$, a crucial element for this particular step.

\begin{lemma}[$L^2$-representation of Fr\'echet derivative]\label{lem:frechet-derivative}
	Let $\bv \in L^2(\Omega; H^1_0(D_T))$ for almost all  $\omega \in \Omega$ be defined as the weak solution to
	\begin{align}\label{def_v_help}
	\left\{
	\begin{array}{rll}
		-\operatorname{div} (\bA \nabla \bv )&=\id_K &\text{in } D_T,\\
		\bv  &= 0 &\text{on } \partial D_T.
	\end{array}\right.
	\end{align}
	The $L^2$-representation of the Fr\'echet derivative of $(\philoc, \id_K)_{L^2(K)}$ is then given by 
	\begin{equation*}
		\frac{\partial}{\partial \bA} (\philoc, \id_K)_{L^2(K)} = - \nabla \philoc \otimes \nabla \bv,
	\end{equation*}
where $\otimes\colon \R^d \times \R^d \rightarrow \R^{d\times d}$ denotes the outer product.
\end{lemma}
\begin{proof}
	Let $\omega \in \Omega$ be fixed. We rewrite the Fr\'echet derivative of $(\philoc, \id_K)_{L^2(K)}$ with respect to $L^2(\mathbb{R}^d;\mathbb{R}^{d\times d})$ as 
	\begin{equation*}
		\label{eq:prooffrechet}
		\frac{\partial}{\partial \bA} (\philoc, \id_K)_{L^2(K)}(\delta \bA) = \left(\frac{\partial \philoc}{\partial \bA}(\delta \bA),\id_K\right)_{L^2(K)}= \int_{D_T} \left(\bA \nabla \bv\right) \cdot \nabla \frac{\partial \philoc}{\partial \bA}(\delta \bA) \dx,
	\end{equation*}
where we tested the  weak formulation of \cref{def_v_help} with $\frac{\partial \philoc}{\partial \bA}(\delta \bA)(\omega)\in H^1_0(D_T)$. 
 	To further simplify the expression on the right-hand side, we differentiate~\cref{eq:basisloc} with respect to $\bA$ using the product rule. This gives for any $w \in H^1_0(D_T)$  that
	\begin{align*}
		\int_{D_T} \big(\delta \bA \nabla \philoc\big)\cdot\nabla w \dx + \int_{D_T} \bA \nabla \frac{\partial \philoc}{\partial \bA}(\delta \bA) \cdot  \nabla w \dx = 0.
	\end{align*}
	Using the test function $w = \bv(\omega) \in H^1_0(D_T)$ and combining the previous two identities, we get 
	\begin{equation*}
		\frac{\partial}{\partial \bA} (\philoc, \id_K)_{L^2(K)}(\delta \bA) = - \int_{D_T}  \big(\delta \bA\nabla \philoc\big)\cdot\nabla \bv \dx.
\end{equation*}
This expression directly characterizes the $L^2$-representation of the Fr\'echet derivative of $(\philoc, \id_K)_{L^2(K)}$, and therefore yields the assertion.
\end{proof}

Another ingredient in the error analysis is the following regularity result for the localized basis functions. The result is needed to further estimate the term we get after applying the spectral gap inequality.  
The proof of this result relies on the condition that the patches take the form of $d$-dimensional bricks, cf.~\cref{lem:RegularityAnnealed}. This condition can be guaranteed, for example, by considering a brick-shaped domain equipped with a Cartesian mesh.

\begin{remark}[Tilde notation]
	In the following, we will write $a \lesssim b$ or $b\gtrsim a$ if it holds that $a \leq C b$ or $a \geq C b$, respectively, where $C>0$ is a constant that may depend on the domain, the shape of the elements and the bounds $\alpha,\beta$ of $\bA$, but is independent of the discretization parameters $H$, $\ell$ and the variations of $\bA$. 
\end{remark}

 \begin{lemma}[$L^4$-regularity of localized basis functions]\label{lem:lemmaregphiloc}
 	Let $\Ab$ be a random coefficient field subject to \cref{assumption_coefficient}. Then, assuming that the patches $D_T$ take the form of bricks, the localized basis functions $\philoc$ satisfy that  	
 	\begin{equation*}
 		\int_{D_T} \mathbb{E}\bigg[\Big(\fint_{B_\varepsilon(x)}  |\nabla \philoc |^2 \dtx \Big)^{2} \bigg] \dx \lesssim \left(\ell H\right)^{4-d}.
 	\end{equation*}
 \end{lemma}

\begin{proof}
In order to apply \cref{lem:RegularityAnnealed}, we need to construct a function $b_T$ such that the localized basis function $\philoc$ is the weak solution to
	\begin{equation*}
		-\nabla\cdot (\Ab\nabla \philoc) = \nabla \cdot b_T \quad \text{on} \;\; D_T
	\end{equation*}
subject to homogeneous Dirichlet boundary conditions on $\partial D_T$.	To this end, one may choose $b_T \coloneqq \nabla r$ for $r$ solving the Laplace problem $\Delta r = g_T$ subject to  homogeneous Dirichlet boundary conditions on $\partial D_T$. 
	With \cref{lem:RegularityAnnealed} we then obtain that
	\begin{equation*}
		\int_{D_T} \mathbb{E}\bigg[\Big(\fint_{B_\varepsilon(x)}  |\nabla \philoc |^2 \dtx \Big)^{2} \bigg] \dx \lesssim \vert D_T\vert^{1-4/q}\Big(\int_{D_T}\vert b_T\vert^q \dx \Big)^{4/q}
	\end{equation*}
	for any $4<q<\infty$. Using standard elliptic regularity on convex domains yields that
	\begin{equation*}
		\|\nabla b_T\|_{L^2(D_T)} = \|D^2 r \|_{L^2(D_T)}\lesssim \| g_T \|_{L^2(D_T)}=1,
	\end{equation*}
	since $g_T$ is $L^2$-normalized. Using the Cauchy--Schwarz inequality and Friedrichs' inequality on $D_T$ for $r\in H^1_0(D_T)$, we get that
	\begin{equation*}
		\|\nabla r\|_{L^2(D_T)}^2 = (g_T,r)_{L^2(D_T)}\leq \|g_T\|_{L^2(D_T)}\|r\|_{L^2(D_T)}\lesssim \ell H  \|g_T\|_{L^2(D_T)}\|\nabla r\|_{L^2(D_T)}.
\end{equation*} With the definition of $b_T$ it follows directly that
	\begin{equation*}
		\|b_T\|_{L^2(D_T)}\lesssim \ell H\|g_T\|_{L^2(D_T)}  = \ell H.
	\end{equation*} 
Applying the Sobolev embedding ($q=6$ is the critical exponent for $d=3$) and a scaling argument (the embedding constant scales with the diameter of~$D_T$), we obtain that 
	\begin{equation*}
		\int_{D_T} \vert b_T \vert^q \dx \lesssim (\ell H)^{d-qd/2}\|b_T\|_{L^2(D_T)}^q + 
		(\ell H)^{d+q(2-d)/2}\|\nabla b_T \|_{L^2(D_T)}^q \lesssim (\ell H)^{d+q(2-d)/2}.
	\end{equation*}
Combining the previous inequalities and setting $q = 5$ gives the assertion.
\end{proof}

The following theorem encapsulates the main result of this work, giving an a posteriori error bound for the proposed numerical stochastic homogenization method. 
\begin{theorem}[A posteriori error bound]\label{t:error_analysis}
	Let $\bA$ be a random coefficient field subject to \cref{assumption_coefficient}. Then, if \cref{riesz_stability} is satisfied, the solution~\cref{eq:stochhommethod} of the proposed numerical stochastic homogenization method satisfies for any $f \in L^2(D)$ that
	\begin{equation*}
		\norml2l2{\ub- \usbar}{D}\lesssim H + C_\mathrm{rb}^{1/2}(H,\ell) \ell^{d/2} \big(\sigma(H,\varepsilon,\ell)  + \varepsilon^{d/2} \ell^2 H^{(4-d)/2}\big) \norm{f}{L^2(D)}
	\end{equation*}
with $C_\mathrm{rb}$ from \cref{riesz_stability}.
\end{theorem}

\begin{proof}
	For the error analysis, we introduce the function 
	\begin{equation}
		\label{eq:dechelpfun}
	\ub_{H,\ell}\coloneqq \sum_{T\in\TH}c_T \philoc,	
	\end{equation}
	where  $(c_T)_{T \in \TH}$ are the coefficients of the representation of $\PiH f$ in terms of the local source terms $\{g_T\with T \in \TH\}$. Using the triangle inequality, we  obtain that
	\begin{align*}
		\begin{split}
			&	\norml2l2{\ub- \usbar}{D} \\
			& \qquad \leq 	\norml2l2{\ub- \PiH\ub}{D} + 	\norml2l2{\PiH(\ub -  \uH)}{D} + 	\norml2l2{\PiH\uH - \usbar}{D} \\
			&\qquad \eqqcolon
			 \Xi_1+\Xi_2+\Xi_3.
		\end{split}
	\end{align*}
In the subsequent analysis, we will estimate the terms $\Xi_1$, $\Xi_2$, and $\Xi_3$ separately. Prior to this, we mention the following approximation result for $\PiH$, the $L^2$-orthogonal projection onto $\TH$-piecewise constants: It holds that
\begin{equation}
	\Vert v- \Pi_H v \Vert_{L^2(T)} \lesssim H  \Vert \nabla v \Vert_{L^2(T)}, \quad v\in H^1(T),\; T \in \TT_H\label{L2-approx-properties};
\end{equation}
see, e.g., \cite{Poincareoriginal, Poincare}.
For the term $\Xi_1$, we obtain using the approximation result~\cref{L2-approx-properties} and the stability estimate~\cref{eq:stabilitysol} that
	\begin{equation*}
		\Xi_1^2 = \E\big[\|\ub- \PiH\ub\|_{L^2(D)}^2\big] \lesssim H^2 \E\big[\|\nabla \ub\|_{L^2(D)}^2\big]\lesssim H^2 \normL2{f}{D}^2.
	\end{equation*}

For estimating the term $\Xi_2$, we first apply the $L^2$-stability of $\Pi_H$ and Friedrichs' inequality. Then,  following the lines of the convergence proof of the SLOD in the deterministic setting, cf.~\cite[Thm.~6.1]{HaPe21b}, we obtain that
	\begin{align*}
		\Xi_2 &
		 \leq C_\mathrm{F}\norml2l2{\nabla (\ub - \uH)}{D} \lesssim (H + C_\mathrm{rb}^{1/2}(H,\ell) \ell^{d/2} \sigma(H,\varepsilon,\ell)) \norm{f}{L^2(D)}.
	\end{align*}

In order to estimate term $\Xi_3$, we recall definitions \cref{eq:dechelpfun,eq:stochhommethod} and use the Cauchy--Schwarz inequality to obtain that
	\begin{align}\label{eq:xi_3}
		\begin{split}
		\Xi_3^2  & = \sum_{T\in\TH} c_T \;\E \Big[\big(\PiH\philoc-\PiH\E[\philoc], \PiH\uH - \usbar\big)_{L^2(D_T)}\Big] \\
		&\leq \sum_{T\in\TH} |c_T|\norml2l2{\PiH\philoc-\PiH\E[\philoc]}{D_T}	\norml2l2{\PiH\uH - \usbar}{D_T}.
		\end{split}
	\end{align} 
Algebraic manipulations then yield for the first term of each summand on the right-hand side of the previous inequality that
	\begin{align}\label{eq:referencephi}
		\begin{split}
		&\norml2l2{\PiH\philoc - \PiH \E [\philoc]}{D_T}^2\\
		&\qquad = \E\bigg[\int_{D_T} \Big(\sum_{K \subset D_T} \big((\philoc,\id_K)_{L^2(K)} - \E[(\philoc, \id_K)_{L^2(K)}]\big)\vert K\vert^{-1} \id_K \Big)^2\dx\bigg] \\
		&\qquad  = \sum_{K \subset D_T}\vert K\vert^{-1}\; \E\Big[\big((\philoc,\id_K)_{L^2(K)} - \E[(\philoc, \id_K)_{L^2(K)}]\big)^2\Big].
			\end{split}
	\end{align}	
	Applying the spectral gap inequality \cref{SpectralGap} and using the $L^2$-representation of the Fr\'echet derivative from \cref{lem:frechet-derivative}, we obtain that
	\begin{align*}
		&\E\Big[\big((\philoc,\id_K)_{L^2(K)} - \E[(\philoc, \id_K)_{L^2(K)}]\big)^2\Big] \\[.5ex]
		&\qquad  \lesssim  \varepsilon^d \; \E\bigg[\int_{D_T}\Big(\fint_{B_\varepsilon(x)} \big\vert \nabla\philoc\otimes\nabla \bv \big\vert \text{d}\tilde{x} \Big)^2\dx\bigg]\\
		&\qquad \leq \varepsilon^d \bigg(\int_{D_T}\E \bigg[\Big(\fint_{B_\varepsilon(x)}\vert \nabla\philoc\vert^2\text{d}\tilde{x}\Big)^2\bigg]\dx\bigg)^{1/2}
		\bigg(\int_{D_T}\E \bigg[\Big(\fint_{B_\varepsilon(x)}\vert \nabla \bv\vert^2\text{d}\tilde{x}\Big)^2\bigg]\dx\bigg)^{1/2},
	\end{align*}
where we used the Cauchy--Schwarz inequality. \cref{lem:lemmaregphiloc} can be employed to bound the first factor on the right-hand side of the preceding inequality. For estimating the second factor, we note that problem \cref{def_v_help} for $\bv\in L^2(\Omega;H^1_0(D_T))$ has the same structure as problem \cref{eq:basisloc} for the localized basis functions. Consequently, a result analogous to \cref{lem:lemmaregphiloc} also holds for $\bv$, leading to
	\begin{align*}
		\int_{D_T}\E \bigg[\Big(\fint_{B_\varepsilon(x)}\vert \nabla \bv\vert^2\text{d}\tilde{x}\Big)^2\bigg]\dx 
		\lesssim (\ell H)^{4-d} \|\id_K\|_{L^2(D_T)}^4.
	\end{align*}
Inserting the estimates for $\philoc$ and $\vb$, we get  that
	\begin{align*}
		\E\Big[\big((\philoc,\id_K)_{L^2(K)} - \E[(\philoc, \id_K)_{L^2(K)}]\big)^2\Big] \lesssim \varepsilon^d (\ell H)^{4-d} \normL2{\id_K}{D_T}^2 = \varepsilon^d (\ell H)^{4-d} \vert K\vert .
	\end{align*}
Using this, we continue to estimate \cref{eq:referencephi} as follows  
	\begin{equation*}
		\label{eq:estpihphiloc}
		\norml2l2{\PiH\philoc - \PiH \E [\philoc]}{D_T}^2
		 \lesssim \varepsilon^d \ell^4 H^{4-d}.
	\end{equation*}
Inserting this estimate into~\cref{eq:xi_3}, applying the Cauchy-Schwarz inequality, recalling the finite overlap of the patches, and utilizing \cref{riesz_stability}, we  finally obtain for $\Xi_3$ that
	\begin{align*}
		\Xi_3^2  & \lesssim  \varepsilon^{d/2} \ell^2 H^{(4-d)/2} \sqrt{\sum_{T\in\TH} c_T^2} \sqrt{\sum_{T\in\TH}\norml2l2{\PiH\uH - \usbar}{D_T}^2}\\
		& \lesssim \varepsilon^{d/2} \ell^{2+d/2} H^{(4-d)/2} C_\mathrm{rb}^{1/2}(H, \ell) \,\Xi_3.
	\end{align*}
The assertion follows immediately after combining the estimates for $\Xi_1$, $\Xi_2$, and $\Xi_3$.
\end{proof}

\section{Error analysis using LOD techniques}
\label{sec:erroranalysiswithlod}
This section utilizes LOD theory to derive an upper bound for the quantity $\sigma$ that appears in the error estimate from \cref{t:error_analysis}. We further estimate $C_\mathrm{rb}$ for the choice of LOD basis functions made for the upper bound on $\sigma$. 

\subsection{Localization error indicator}

We first derive an upper bound for the localization error $\sigma$ defined in~\cref{eq:locerrind}. The bound is based on the lowest-order LOD from~\cite{Mai20ppt,HaPe21,DHM23}, whose construction uses non-negative bubbles $\{b_T \with T \in \TH\}$; see also \cite{FeP20}. The bubble function $b_T\in H_0^1(T)$ is chosen such that $\Pi_Hb_T = \id_T$ and
\begin{equation}
	\label{eq:bubble}
	\Vert b_T \Vert_{L^2(T)} \lesssim H \Vert \nabla b_T \Vert_{L^2(T)} \lesssim \sqrt{\vert T \vert}
\end{equation}
holds. Recalling the abbreviation $D_T = \mathsf N^\ell(T)$ for the $\ell$-th order patch around $T$, cf.~\cref{eq:patch}, we introduce the space of fine-scale functions supported on~$D_T$ by $\mathcal{W}_{T,\ell} \coloneqq \{w\in H_0^1(D_T) \with \Pi_{H, D_T} w = 0\}$. 
The LOD basis function corresponding to the element $ T\in\TH$ is then defined by 
\begin{equation}
	\label{eq:lodbasis}
	\philod \coloneqq (1-\corrc)b_T \in L^2(\Omega; H^1_0(D_T)),
\end{equation}
where $\corrc b_T \in L^2(\Omega;\mathcal W_{T,\ell})$ denotes the fine-scale correction of the bubble $b_T$, which is defined for almost all $\omega \in \Omega$~by 
\begin{equation}\label{eq:corrc}
	\ba(\corrc b_T, w) = \ba(b_T,w) \quad \text{ for all } w\in\mathcal{W}_{T,\ell}.
\end{equation}
Note that the well-posedness of the operator $\corrc$ is a consequence of the Lax--Milgram theorem, recalling that $\mathcal{W}_{T,\ell}$ is a closed subspace of $H^1_0(D_T)$.

In the following lemma we derive an upper bound on $\sigma$, based on the observation that the LOD basis function~$\philod$ possesses a $\mathcal T_H$-piecewise constant source term
\begin{equation*}
	\glod \coloneqq - \mathrm{div} \Ab \nabla \philod \in L^2(\Omega;\Pnull(\T_{H,D_T}));
\end{equation*}
see, e.g.,~\cite[Lem.~A.2]{HaPe21b}.

\begin{lemma}[Upper bound on $\sigma$]\label{lemma:ub}
	Choosing an $L^2$-normalized version of $g_T  \coloneqq \E[\glod]$ in 
	\cref{eq: sigma} yields the upper bound 
		\begin{equation}
			\label{eq:estsigma}
		\sigma \lesssim  \ell^2H^{-1} \exp(-C_\mathrm{d} \ell) + \ell^4\left(\frac{\varepsilon}{H}\right)^{d/2}
	\end{equation}
	with $C_\mathrm{d}>0$ independent of $H$ and $\ell$, provided that $\varepsilon$ satisfies the smallness assumption
	\begin{equation}\label{assumption_eps1}
		\varepsilon^d \lesssim \ell^{-8}H^d.
	\end{equation}
\end{lemma}

\begin{proof}
For all  $\bv \in \Yb \subset L^2(\Omega; H^1_\Gamma(D_T))$ it holds that  $\E[\trinv \tr \bv] = \E[\bv]$. Hence, by inserting $g_T = \E[\glod]$ into \cref{eq: sigma}, we obtain that
	\begin{equation*}
		\sigma_T(H,\varepsilon,\ell) \leq \frac{1}{\| g_T \|_{L^2(D_T)}}\; \underset{\substack{\bv\in L^2(\Omega;H^1_\Gamma(D_T)) \\ \Vert \bv \Vert_{L^2(\Omega;H^1(D_T))}=1}}{\sup} (g_T,\E[\trinv \tr \bv])_{L^2(D_T)}.
	\end{equation*}
Note that by dividing by the norm of $g_T$, we account for the fact that $g_T$ may not be normalized. We denote by $\Apatchinv\colon L^2(\Omega;L^2(D_T))\to L^2(\Omega;H^1_0(D_T))$ the local solution operator defined on the patch $D_T$, which satisfies the 
following stability estimate
\begin{equation}
	\label{eq:stabilitysolop}
	\norml2l2{\nabla \Apatchinv \bg}{D_T} \lesssim \norml2l2{\bg}{D_T}.
\end{equation}
Therefore, we obtain for any $\bv\in L^2(\Omega; H^1_\Gamma(D_T))$ that
	\begin{align*}
		(g_T,\E[\trinv \tr \bv])_{L^2(D_T)} &=  \E[(\glod, \bv)_{L^2(D_T)} - \ba(\Apatchinv\glod,\bv)]  \\&\qquad + \E[(g_T - \glod, \bv)_{L^2(D_T)} - \ba(\Apatchinv (g_T-\glod),\bv)]\\
		&\eqqcolon \Xi_1 + \Xi_2.
	\end{align*}

To estimate the term $\Xi_1$, we apply the deterministic result  \cite[Lem.~6.4]{HaPe21b} for any $\omega \in \Omega $ and use the Cauchy--Schwarz inequality to get that
	\begin{equation*}
		\Xi_1\lesssim H^{-1} \exp(-C_\mathrm{d} \ell) \|\bv\|_{L^2(\Omega,H^1(D_T))}\,\norml2l2{\glod}{D_T},
	\end{equation*}
where $C_\mathrm{d}>0$ is independent of $H$ and $\ell$. Using the estimate
	\begin{equation}\label{eq:norm_glod}
		\norml2l2{\glod}{D_T} \lesssim H^{d/2-2},
	\end{equation}
which can be derived by taking the expectation of the corresponding deterministic identity from \cite[Lem.~A.2]{HaPe21b}, yields that 
	\begin{align*}
		\Xi_1 \lesssim  H^{d/2-3} \exp(-C_\mathrm{d} \ell)\, \|\bv\|_{L^2(\Omega,H^1(D_T))}.
	\end{align*}

For the term $\Xi_2$, we obtain using \cref{eq:stabilitysolop}   
and the Cauchy--Schwarz inequality that 
	\begin{align*}
		\Xi_2 \lesssim \norml2l2{g_T - \glod}{D_T} \|\bv\|_{L^2(\Omega,H^1(D_T))}.
	\end{align*}
In order to estimate the first factor on the right-hand side, we proceed similarly as in the proof of \cref{t:error_analysis} to obtain that
	\begin{align*}
		 \norml2l2{g_T - \glod}{D_T}^2 
		= \sum_{K\subset D_T}\vert K\vert^{-1}\; \E\Big[\big((\glod, \id_K)_{L^2(K)} - \E\left[(\glod, \id_K)_{L^2(K)}\right]\big)^2\Big].
	\end{align*}	
	Using the spectral gap inequality \cref{SpectralGap}, we obtain for each summand that
	\begin{align}
		\label{eq:rhsspectralgapie}
		\begin{split}
			&\E\Big[\big((\glod,\id_K)_{L^2(K)} - \E[(\glod, \id_K)_{L^2(K)}]\big)^2\Big]\\*
			& \qquad \lesssim \varepsilon^d \; \E\bigg[\int_{\R^d} \Big(\fint_{B_\varepsilon(x)} \Big|\frac{\partial (\glod, \id_K)_{L^2(K)}}{\partial \Ab}(\tilde{x})\Big| \mathrm{d}\tilde{x} \Big)^2 \dx\bigg].
		\end{split}
	\end{align}
	The $L^2$-representation of the Fr\'echet derivative of  $(\glod,\id_K)_{L^2(K)}$ is derived in \cref{lem:frech-glod}. It consists of a sum  of outer products of the gradients of combinations of $b_T$, $b_K$, $\corrc b_T$ and $\corrc b_K$. To estimate the summands involving bubble functions, we utilize the property~\cref{eq:bubble} for all $K \subset D_T$ and derive the estimate
	\begin{align}\label{eq:bubbleest}
	\int_{D_T} \Big(\fint_{B_\varepsilon(x)} |\nabla  b_K|^2 \dtx\Big)^2 \dx  \lesssim H^{d-4}.
	\end{align}
	To proceed with the estimation of \cref{eq:rhsspectralgapie}, we need to estimate the four terms resulting from the summands of the Fr\'echet derivative, cf.~\cref{lem:frech-glod}. In the following, we present the estimate for the second term, noting that all other estimates follow analogously. By employing the regularity result from \cref{l:lemmaregc} and \cref{eq:bubbleest}, we obtain that
	\begin{align*}
		&\mathbb E \bigg[\int_{D_T} \Big(\fint_{B_\varepsilon(x)}| \nabla \corrc b_T \otimes \nabla b_K| \dtx\Big)^2 \dx\bigg]\\
		& \quad \leq  \bigg(\int_{D_T}\mathbb E \bigg[ \Big(\fint_{B_\varepsilon(x)} |\nabla \corrc b_T|^2 \dtx\Big)^2 \bigg] \dx \bigg)^{1/2}
		\bigg(\int_{D_T} \Big(\fint_{B_\varepsilon(x)}|\nabla b_K |^2 \dtx\Big)^2 \dx\bigg)^{1/2}\\
		&\quad \lesssim  \ell^{2-d/2}H^{d-4},
	\end{align*}
	where we used the Cauchy--Schwarz inequality. Note that all four terms can be majorized by $\ell^{4-d}H^{d-4}$, which results from estimating the last summand. The combination of the previous estimates yields that
	\begin{align}\label{eq: diff_g_glod}
	 \norml2l2{g_T - \glod}{D_T} \lesssim \bigg(\sum_{K\subset D_T }\vert K\vert^{-1}\; \varepsilon^d \ell^{4-d}H^{d-4}\bigg)^{1/2} \lesssim \varepsilon^{d/2}\ell^2H^{-2}.
	\end{align}

Using the estimate 
\begin{equation*}
	\norml2l2{\glod}{D_T} \gtrsim \ell^{-2}H^{d/2-2},
\end{equation*}
which can be derived by taking the expectation 
of the corresponding deterministic identity from \cite[Lem.~A.2]{HaPe21b}, we can derive the following lower bound for the $L^2$-norm of $g_T$
\begin{align}\label{eq: norm_g}
	\begin{split}
		\| g_T \|_{L^2(D_T)}^2 &= \norml2l2{ g_T}{D_T}^2 \\
		&\geq \frac{1}{2}\norml2l2{\glod }{D_T}^2  - \norml2l2{g_T-\glod }{D_T}^2 \\
		& \gtrsim \frac{1}{2}\,\ell^{-4}H^{d-4} - \ell^4H^{-4} \varepsilon^{d} \gtrsim  \ell^{-4}H^{d-4}.
	\end{split}
\end{align}
Here, we used the reverse triangle inequality, the weighted Young's inequality for showing that for $a,b\geq0$ it holds that $|a-b|^2\geq \tfrac{a^2}{2}-b^2$, as well as the smallness assumption~\cref{assumption_eps1}. Finally, combining all estimates leads to
\begin{align*}
	\sigma_T \lesssim  \frac{1}{\| g_T \|_{L^2(D_T)}}\big(H^{d/2-3} \exp(-C_\mathrm{d} \ell)+ \ell^2H^{-2} \varepsilon^{d/2}\big)\lesssim \ell^{2}H^{-1} \exp(-C \ell) + \ell^4\left(\frac{\varepsilon}{H}\right)^{d/2}.
\end{align*}
The assertion follows directly when taking the maximum over all $T \in \TH$.
\end{proof}

Combining this a priori result for $\sigma$ with \cref{t:error_analysis} yields the error estimate given in the following corollary. The Riesz constant $\Crb$ can be computed a posteriori, cf.~\cref{sec:implementation}.
\begin{corollary}[Combined error bound]\label{c:error_analysis_combined}
	Suppose that the assumptions of \cref{,t:error_analysis,lemma:ub} are fulfilled and that $\ell\gtrsim |\log H|$ holds. Then, 
	the solution~\cref{eq:stochhommethod} of the proposed numerical stochastic homogenization method satisfies, for any $f \in L^2(D)$, that
	\begin{equation*}
		\norml2l2{\ub- \usbar}{D}\lesssim H + C_\mathrm{rb}^{1/2}(H,\ell) \ell^{4+d/2} \left(\frac{\varepsilon}{H}\right)^{d/2}\norm{f}{L^2(D)}.
	\end{equation*}
\end{corollary}

\subsection{Riesz stability}
In a next step, we show that the local source terms corresponding to the LOD basis functions \cref{eq:lodbasis} are Riesz stable in the  sense of \cref{riesz_stability}.

\begin{lemma}[Riesz stability of LOD source terms]\label{lem:rieszlod}
	Suppose that $\ell$ is chosen such that $\ell \gtrsim |\log(H)|$ and that $\varepsilon$ satisfies the smallness assumption
	\begin{equation}\label{assumption_eps2}
		\varepsilon^d \lesssim\ell^{-(8+d)}H^{4+d}.
	\end{equation} 
Then, for the local source terms $g_T = \E[\glod]$ it holds for all $(c_T)_{T\in\TH}$ that
	\begin{equation}
		\label{eq:reisz}
		H^4 \sum_{T\in\TH} c_T^2 \lesssim \bigg\|\sum_{T\in\TH} c_T \frac{g_T}{\norm{g_T}{L^2(D_T)}}\bigg\|_{L^2(D)}^2.
	\end{equation}
\end{lemma} 
\begin{proof}
	We begin the proof by noting that applying the weighted Young inequality twice gives the elementary estimate $|a-b-c|^2 \geq \tfrac14|a|^2-|b|^2-|c|^2$ for any $a,b,c\geq0$. Combining this with the inverse triangle inequality, we obtain that
	\begin{align*}
		&\bigg\|\sum_{T\in\TH} c_T \frac{g_T}{\norm{g_T}{L^2(D_T)}}\bigg\|_{L^2(D)}^2 
		=\bigg\|\sum_{T\in\TH} c_T \frac{g_T}{\norm{g_T}{L^2(D_T)}}\bigg\|_{L^2(\Omega,L^2(D))}^2 \\
		&\qquad  \geq \frac{1}{4}\,\bigg\|\sum_{T\in\TH} c_T \frac{\glod}{\norm{\glod}{L^2(D_T)}}\bigg\|_{L^2(\Omega,L^2(D))}^2 - \bigg\|\sum_{T\in\TH} c_T \frac{g_T - \glod}{\norm{g_T}{L^2(D_T)}}\bigg\|_{L^2(\Omega,L^2(D))}^2\\*
		& \qquad\qquad  - \bigg\|\sum_{T\in\TH} c_T \bigg(\frac{\glod}{\norm{g_T}{L^2(D_T)}} - \frac{\glod}{\norm{\glod}{L^2(D_T)}}\bigg)\bigg\|_{L^2(\Omega,L^2(D))}^2 \\*
		& \qquad \eqqcolon \frac{1}{4}\,\Xi_1 - \Xi_2 - \Xi_3. 
	\end{align*}
For estimating the term $\Xi_1$ from below, we use the corresponding deterministic result from \cite[Lem.~6.4]{HaPe21b} and take the expectation which yields that
	\begin{equation*}
		\Xi_1 \gtrsim H^4 \sum_{T\in\TH} c_T^2.
	\end{equation*} 
To estimate the term $\Xi_2$ from above, we use the finite overlap of the patches $D_T$ as well as estimates \cref{eq: diff_g_glod,eq: norm_g} to get that
	\begin{align*}
		\Xi_2 &\lesssim \ell^{4+d} H^{4-d} \sum_{T\in\TH} c_T^2 \, \norml2l2{g_T - \glod}{D_T}^2 \lesssim \ell^{8+d} \varepsilon^dH^{-d} \sum_{T\in\TH} c_T^2.
	\end{align*}
The estimate for $\Xi_3$ can be derived similarly using again the finite overlap of the patches~$D_T$, the reverse triangle inequality, \cref{eq: diff_g_glod,eq: norm_g}. We obtain that
	\begin{align*}
		\Xi_3 &=\bigg\|\sum_{T\in\TH} c_T \frac{\glod(\norm{\glod}{L^2(D_T)}- \norm{g_T}{L^2(D_T)})}{\norm{g_T}{L^2(D_T)}\norm{\glod}{L^2(D_T)}} \bigg\|_{L^2(\Omega,L^2(D))}^2 \\
		& \lesssim \ell^d \sum_{T\in\TH} c_T^2 \, \E\Bigg[ \frac{\norm{g_T - \glod}{L^2(D_T)}^2}{\norm{g_T}{L^2(D_T)}^2}  \Bigg] \lesssim \ell^{8+d} \varepsilon^dH^{-d}  \sum_{T\in\TH} c_T^2.
	\end{align*}
Combining the previous estimates and using the smallness assumption \eqref{assumption_eps2} yields the assertion.
\end{proof}

\section{Practical implementation}\label{sec:implementation}

To effectively implement the proposed numerical stochastic homogenization method, it is crucial to employ an efficient sampling strategy for the space $\Yb$ and ensure that the local source terms $\{g_T\with T\in \TH\}$ form a stable basis of $\mathbb{P}^0(\TH)$. These aspects will be addressed in the following two subsections.

\subsection{Sampling of the space $\Yb$}
We consider an arbitrary patch $D_T$ and denote the number of coarse elements in this patch by $N \coloneqq \#\TT_{H, D_T}$. In a practical implementation, all local infinite-dimensional problems that appear in the derivation of the basis functions must be replaced by finite-dimensional counterparts. To obtain these finite-dimensional counterparts, we perform a discretization using the $\mathcal Q^1$-finite element method with respect to the fine mesh~$\mathcal{T}_{h,D_T}$ constructed by uniform refinements of~$\TT_{H,D_T}$. The number of elements of $\mathcal{T}_{h,D_T}$ is denoted by $n$. 

To handle the stochasticity in the definition of~$\Yb$, our implementation draws $M$ samples of the random coefficient $\Ab$ and, for each sample, closely follows the methodology outlined in \cite[App.~B]{HaPe21b} for the deterministic case. Specifically, we generate a matrix $\mathbf S_i \in \mathbb{R}^{n\times m}$ for $i=1,\dots,M$, whose columns represent the coordinate vectors of the discrete $\Ab(\omega_i)$-harmonic extensions of $m\in\mathbb N$ samples of random boundary data on $\partial D_T \backslash \partial D$. Then we compute the matrices $\mathbf P_i\in\mathbb R^{N\times m}$ by applying the $L^2$-orthogonal projection onto the characteristic functions $\lbrace \id_T\with T\in\mathcal T_{H,D_T} \rbrace$ column by column to $\mathbf S_i$. Finally, the SVD of the matrix $\mathbf X \coloneqq [\mathbf P_1,\dots,\mathbf P_M]$ is computed, yielding coordinate vectors of potential right-hand sides $g_T$. For details on the practical realization of this SVD, we refer to~\cite[App.~B]{HaPe21b}. Finally, the localized deterministic basis functions are computed as empirical means, again using $M$ samples of the random coefficient.
In the numerical experiments performed in \cref{sec:numexp}, the number of random boundary samples is set to $m=3N$. For the number of random coefficient samples, we use $M = 5000$.

\subsection{Stable local source terms}
Next, we discuss how the stability of the local source terms $\{g_T \with T \in \mathcal T_H\}$ can be ensured in a practical implementation. Our implementation achieves stability by an additional optimization step, similar to the one used in~\cite{Bonizzoni-Freese-Peterseim}. Given the singular values $\sigma_1\geq\sigma_2\geq\dots \geq \sigma_N\geq 0$ of the matrix $\mathbf X$ associated with the patch $D_T$, we consider all indices $1\leq i \leq N$ such that
\begin{align*}\label{eq:choicelocrhs}
	\frac{\sigma_i}{\sigma_1}\leq \max \Big\{\Big(\frac{\sigma_N}{\sigma_1}\Big)^{1/p},10^{-10}\Big\}
\end{align*}
and denote the resulting set of indices by $\mathcal I$. Each index in the set $\mathcal I$ corresponds to a potential candidate for a local source term. For the choice $p = 1$ only the smallest singular value is considered. Since our optimization problem is meaningful whenever multiple functions are considered, we restrict ourselves to the choices $p > 1$.

Among these candidate functions, we choose the one that maximizes a weighted $L^2(D_T)$-norm under the unit mass constraint. The weighted $L^2(D_T)$-norm is defined using a piecewise constant weighting function that is zero in the central element $T$ and grows polynomially as the distance from the center increases. This enforces a concentration of mass in the center of each patch, resulting in linearly independent local source terms $\{g_T \with T \in \TH\}$ in practice. More specifically, we introduce the distance function $\mathrm{dist}(T, K)$ between the elements $T,K \in \mathcal T_H$ as
\begin{align*}
	\mathrm{dist}(T, K) \coloneqq H^{-1}|m_K - m_T | \in \mathbb{N}^d,
\end{align*}
where $m_T , m_K \in \mathbb{R}^d$ are the midpoints of the elements $T$ and $K$, respectively. The weighting function is then defined for each element $K\in\TT_{H, D_T} $ as
\begin{align*}
	w_T(K) \coloneqq \big\vert \mathrm{dist}(T,K)\big\vert_\infty^{r}
\end{align*} 
for a parameter $r\geq1$, where $|\cdot|_\infty$ denotes the infinity norm on $\mathbb R^d$. \cref{fig: weight_function} provides an illustration of this weighting function in two spatial dimensions. In our numerical experiments in \cref{sec:numexp}, we use $p=1.5$ and $r = 6$.

\begin{remark}[Computation of $\Crb$]
	Given the local source terms $\{g_{T_i}\with i=1,\dots,\# \mathcal T_H\}$,  the Riesz stability constant $\Crb$ appearing in \cref{c:error_analysis_combined} equals the reciprocal of the smallest eigenvalue of the matrix $G\in \mathbb{R}^{\#\mathcal T_H \times \#\mathcal T_H}$ with entries given by $G_{ij}=(g_{T_i},g_{T_j})_{L^2(D)}$.
\end{remark}

\begin{figure}
	\begin{center}
		\begin{TAB}(e,0.9cm,0.9cm){|c|c|c|c|c|}{|c|c|c|c|c|}
			64 & 64 & 64 & 64 & 64 \\
			64 & 1 & 1 & 1 & 64 \\
			64 & 1 & 0 & 1 & 64 \\
			64 & 1 & 1 & 1 & 64 \\
			64 & 64 & 64 & 64 & 64   
		\end{TAB}
	\end{center}
\caption{Piecewise constant weighting function $w_T$ for an interior element $T$ with $\ell = 2$ in two spatial dimensions.}
\label{fig: weight_function}
\end{figure}

\begin{remark}[Uniform Cartesian meshes]
Note that in the case of uniform Cartesian meshes, the computational complexity of the method can be significantly reduced when utilizing the stationarity of the coefficient $\Ab$, cf.~\cref{assumption_coefficient}. In fact, only $\mathcal{O} (\ell^d)$ reference patches need to be considered for the computation of the basis functions and local source terms of the method. All other basis functions and local source terms can then be obtained by translation; see, e.g., \cite{GaP15}.
\end{remark}

\section{Numerical experiments}\label{sec:numexp}

The following numerical experiments are intended to demonstrate the effectiveness of the proposed numerical homogenization method. In our implementation, we consider uniform Cartesian meshes of the domain $D = (0,1)^d$ with $d\in\{1,2\}$. Note that from now on we use $H$ to denote the side length of the elements instead of their diameter.
For the solution of the local patch problems and the computation of the reference solution $\ub_h$ we employ the $\mathcal Q_1$-finite element method on the fine mesh $\TT_h$ with $h=2^{-10}$. We denote by $\bar{u}_{H,h,\ell}$ the fully discrete numerical approximation to $\E[\ub]$. In the following all expected values are replaced by appropriate empirical means.

The random coefficients $\Ab$ that are considered in the following numerical experiments are piecewise constant with respect to the uniform Cartesian meshes $\mathcal T_\varepsilon$ with mesh sizes $\varepsilon \in  \lbrace 2^{-5},2^{-6},2^{-7},2^{-8},2^{-9}\rbrace$. These coefficients take independent and identically distributed element values in the interval $[0.1,1]$. We further consider the sequence of coarse meshes~$\TH$ with mesh sizes $H\in \lbrace 2^{-3},2^{-4},2^{-5},2^{-6}\rbrace$. Note that we only consider coarse mesh sizes~$H>\varepsilon$ for which the coarse mesh does not resolve the minimal length scale of the random coefficient. We also exclude combinations of $H$ and $\ell$ for which a patch coincides with the whole domain $D$. To calculate the reference solution, we employ $M=5000$ samples, which is consistent with the number used for the local patch problems. The samples are obtained by a quasi-Monte Carlo sampling strategy in one spatial dimension and a Monte Carlo sampling strategy in two spatial dimensions.

\subsection*{Numerical investigation of $\sigma$ and $C_\mathrm{rb}$}

\begin{figure}
	\includegraphics[width=.49\linewidth]{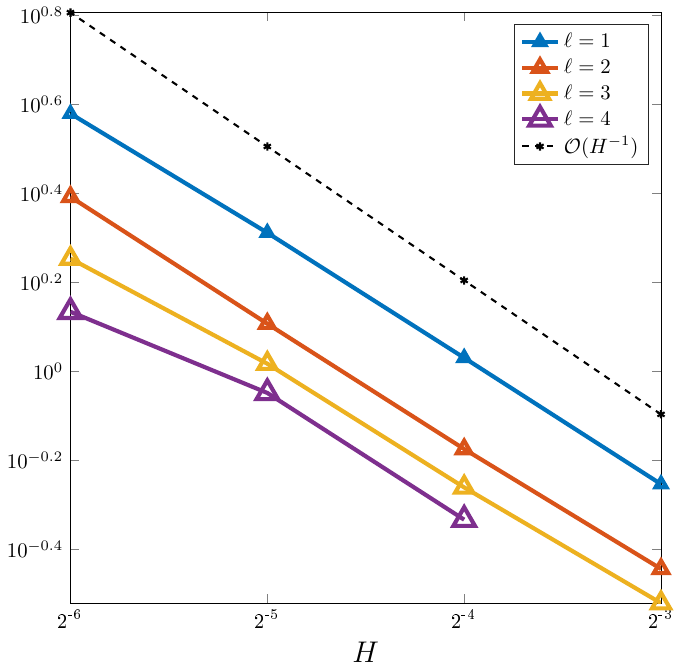}\hfill
	\includegraphics[width=.49\linewidth]{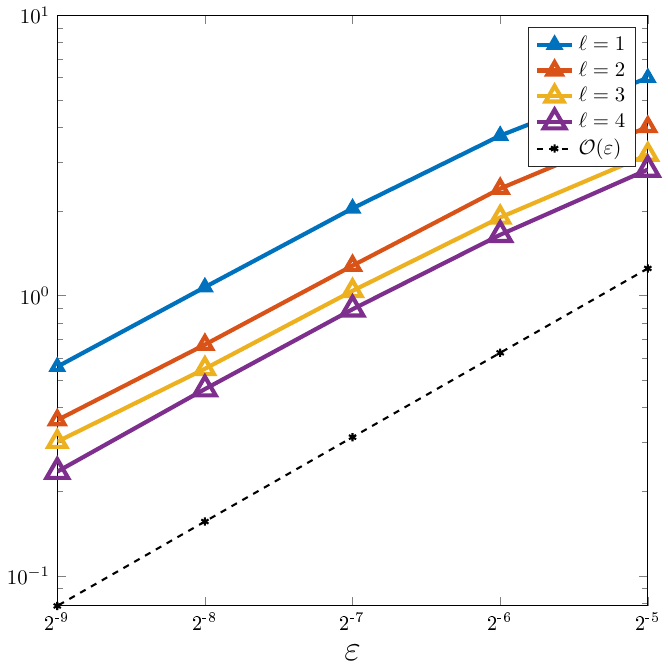}
	\caption{Depiction of $\sigma$ for a $\mathcal T_\varepsilon$-piecewise constant random coefficient in two spatial dimensions. Left: in dependence of the coarse mesh size $H$ for $\varepsilon=2^{-8}$; Right: in dependence of the correlation length $\varepsilon$ for $H=2^{-4}$.}
	\label{fig:sv_d2}
\end{figure}

We first examine the behavior of the localization error indicator $\sigma$ as a function of the coarse mesh size $H$ and the correlation length~$\varepsilon$. For this, we consider the case $d = 2$ and utilize the sequences of coarse meshes and correlation lengths mentioned above. \cref{fig:sv_d2} visualizes the values of~$\sigma$ for a fixed correlation length~$\varepsilon$ and varying mesh sizes $H$ (left) and for fixed $H$ and varying~$\varepsilon$ (right). In both cases one observes a scaling like $\tfrac{\varepsilon}{H}$, which numerically validates the upper bound for $\sigma$ from \cref{lemma:ub} in the case $d = 2$. Note that the stochastic errors dominate, and consequently, the first term in \cref{eq:estsigma}, which decays exponentially in $\ell$, is not visible. Plotting $\sigma$ as a function of $\ell$ would give a scaling like $\ell^{-1/2}$. 

\begin{figure}
	\includegraphics[width=.49\linewidth]{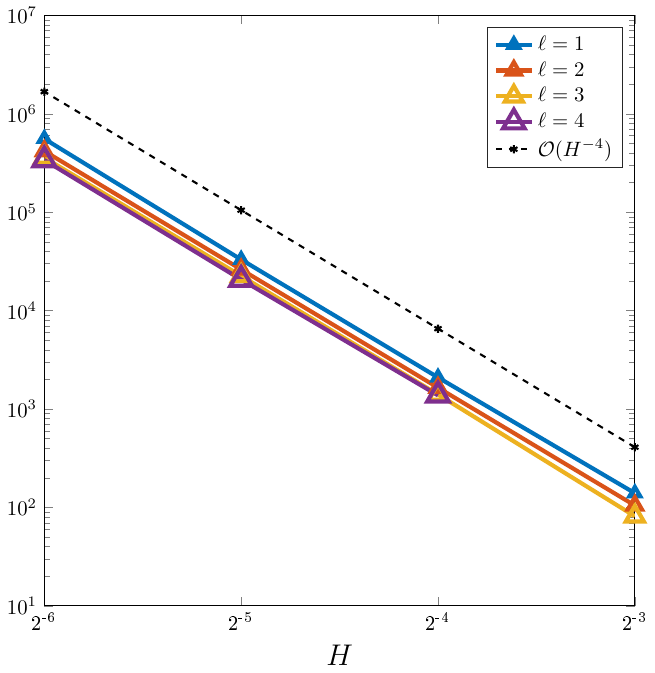}
	\caption{Depiction of the Riesz stability constant $C_\mathrm{rb}$ of the stochastic SLOD as a function of the coarse mesh size $H$ for a $\mathcal T_\varepsilon$-piecewise constant random coefficient with $\varepsilon=2^{-8}$ in two spatial dimensions.}
	\label{fig:riesz_lower}
\end{figure}

Next we examine the behavior of the Riesz stability constant $\Crb$ of the local SLOD source terms as a function of $H$. In \cref{fig:riesz_lower} we observe that $\Crb$ scales like $H^{-4}$, which is consistent with the results for the stochastically averaged LOD source terms proved in \cref{lem:rieszlod}. Our numerical experiments indicate no dependency of the Riesz stability constant on $\varepsilon$ or $\ell$, which is also in line with the findings from \cref{lem:rieszlod}.


\subsection*{Numerical validation of convergence}
\begin{figure}
	\includegraphics[width=.49\linewidth, height=0.475\linewidth]{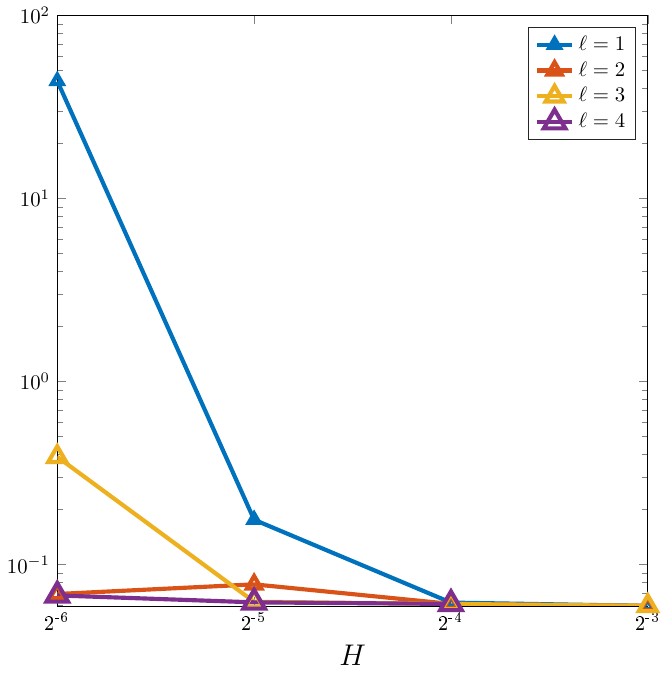}\hfill
	\includegraphics[width=.49\linewidth]{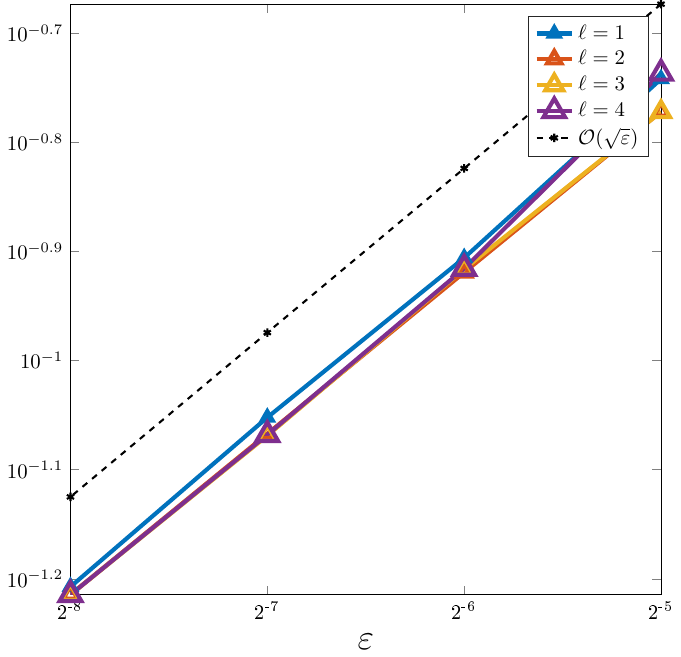}
	\caption{Plot of the relative $L^2$-errors $\|\PiH\ub_h - \bar{u}_{H,h,\ell}\|_{L^2(\Omega;L^2(D))}$ of the proposed SLOD method for a $\mathcal T_\varepsilon$-piecewise constant random coefficient in one spatial dimension. Left: errors as functions of the coarse mesh size $H$ for fixed $\varepsilon = 2^{-8}$ and several oversampling parameters $\ell$; Right: errors in dependency of the correlation length $\varepsilon$ for fixed $H=2^{-4}$ and several values of $\ell$.}
	\label{errL2Hd1_constf}
\end{figure}

\begin{figure}
	\includegraphics[width=.49\linewidth]{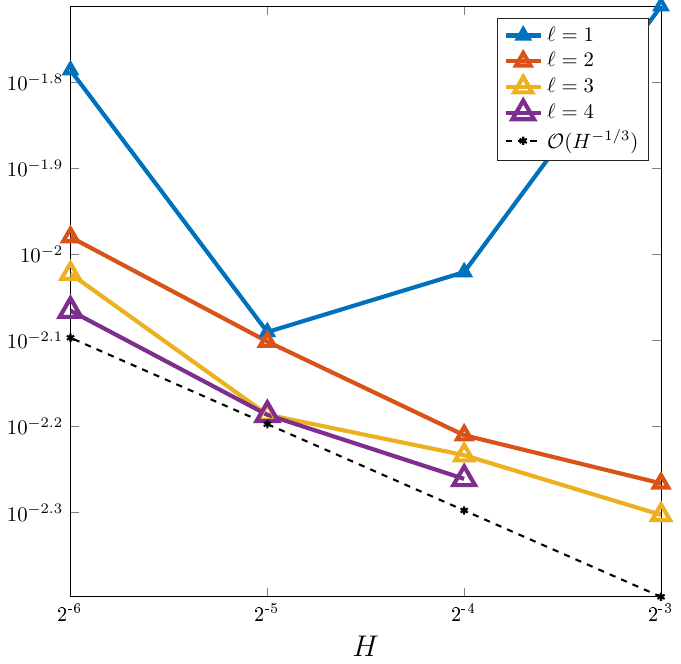}\hfill
	\includegraphics[width=.49\linewidth]{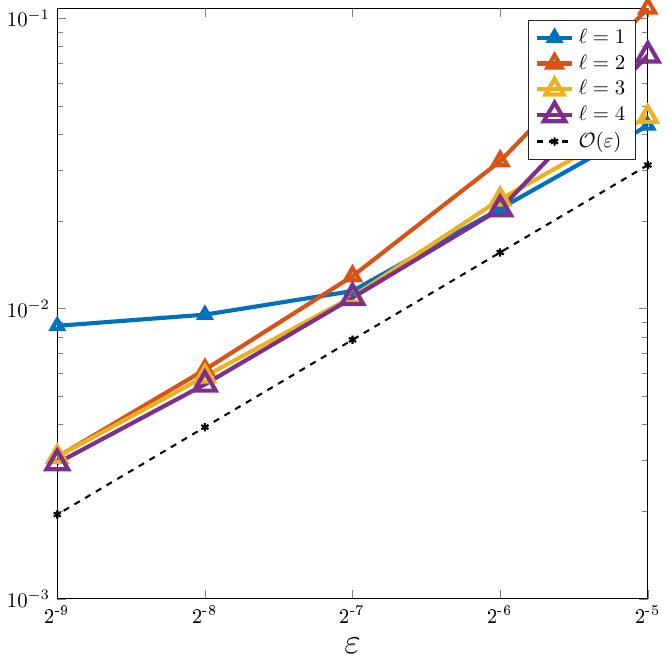}
	\caption{Plot of the relative $L^2$-errors $\|\PiH\ub_h - \bar{u}_{H,h,\ell}\|_{L^2(\Omega;L^2(D))}$ of the proposed SLOD method for a $\mathcal T_\varepsilon$-piecewise constant random coefficient in two spatial dimensions. Left: errors as functions of the coarse mesh size~$H$ for fixed $\varepsilon = 2^{-8}$ and several oversampling parameters $\ell$; Right: errors as functions of the correlation length $\varepsilon$ for fixed $H=2^{-4}$ and several~$\ell$.}
	\label{errL2Hd2_constf}
\end{figure}

To numerically verify the convergence of the proposed numerical stochastic homogenization method, we consider the source terms
\begin{equation*}
	f(x)= 2\pi^2\sin(x),\quad f(x,y)= 2\pi^2\sin(x)\sin(y)
\end{equation*}
in one and two spatial dimensions, respectively. 
\cref{errL2Hd1_constf,errL2Hd2_constf} show the resulting relative $L^2$-errors computed using the reference solution $\ub_h$. For fixed $H$ and varying~$\varepsilon$ we observe the rate $\varepsilon^{d/2}$, which is in agreement with \cref{c:error_analysis_combined}. When considering the converse case, we have to distinguish between one and two spatial dimensions. In one dimension, the expected negative power of $H$ does not manifest itself, and in our numerical experiments the error remains relatively constant with respect to $H$ (provided the coarse mesh is sufficiently coarse compared to $\varepsilon$). In the two-dimensional case, we observe a negative dependence on $H$, which is much weaker than the $H^{-2}$ predicted by \cref{c:error_analysis_combined}. The error rather seems to scale like $H^{-1/3}$.


\appendix
\section{Auxiliary results}

The error analysis of the proposed numerical stochastic homogenization method is based on the so-called Calderon--Zygmund estimates, which are a popular tool in the theory of quantitative stochastic homogenization. Such estimates were established for an equation on the full space $\R^d$ in \cite{Otto2020}, extending earlier results from \cite{ARMSTRONG2016312, DGO2020}.  For annealed Calderon--Zygmund estimates we refer to \cite{JO22,WangXu24}, where the latter work considers the case of Lipschitz domains. 
Contrary to the Calderon--Zygmund estimate given below, these annealed estimates involve only a loss in stochastic integrability and not in spatial integrability. Since such annealed estimates only lead to better (hidden) constants in the final error bounds, we will henceforth stick to a suboptimal Calderon--Zygmund estimate similar to \cite[Lem.~4.8]{FisGP19ppt}, where an a priori error analysis for a related numerical stochastic homogenization method is performed. The proof of the following estimate, which is beyond the scope of this manuscript, is analogous to the full-space case \cite[Thm. 6.1]{Otto2020} and uses the boundary regularity theory of \cite{Fischer2017,JRS24} as well as a classical regularity theory at edges and corners.

\begin{lemma}[Annealed large-scale $L^p$ regularity]
	\label{lem:RegularityAnnealed}
	Let $d \in \lbrace2,3\rbrace$, and let $\Ab$ be a random coefficient field subject to \cref{eq:coeffspd} and \cref{assumption_coefficient}. Let $Q\subset \mathbb{R}^d$ be a box, let $ \bb\in L^2(\Omega; L^2(Q))$, and let $\boldsymbol{u}\in L^2(\Omega;H^1_0(Q))$ be a solution to the linear elliptic PDE
	\begin{equation*}
		\label{e:rhs}
		\begin{split}
			\begin{aligned}
				-\nabla \cdot (\Ab\nabla \boldsymbol{u})&=\nabla \cdot  \bb \quad &&\text{on }Q,\\
				\boldsymbol{u}&\equiv 0 &&\text{on }\partial Q.
			\end{aligned}
		\end{split}
	\end{equation*}
	Then for any $2\leq p<\infty$ and any $p<q<\infty$
	there holds a regularity estimate of the form
	\begin{align*}
		\fint_Q \mathbb{E}\Bigg[\bigg(\fint_{B_\varepsilon(x)}  |\nabla \boldsymbol{u}|^2 \,d\tilde x\bigg)^{p/2} \Bigg] \dx
		\leq C(\lambda,\Lambda,\rho,p,q) \bigg(\fint_Q \mathbb{E}\Bigg[\bigg(\fint_{B_\varepsilon(x)}  |\boldsymbol{b}|^2 \,d\tilde x\bigg)^{q/2} \Bigg] \dx \bigg)^{p/q}.
	\end{align*}
\end{lemma}

In the following, we present two results used in the proof of \cref{lemma:ub}. The first result provides a $L^2$-representation of the Fr\'echet derivative, which is needed to apply the spectral gap inequality.
\begin{lemma}[Fr\'echet derivative of LOD right-hand sides]\label{lem:frech-glod}
	The $L^2$-representation of the Fr\'echet derivative of $(\glod, \id_K )_{L^2(K)}$ is given by 
	\begin{equation*}
		\frac{\partial}{\partial \Ab} (\glod, \id_K )_{L^2(K)}  = \nabla b_T \otimes \nabla b_K - \nabla \corrc b_T \otimes \nabla b_K - \nabla b_T \otimes \nabla \corrc b_K + \nabla \corrc b_T \otimes \nabla \corrc b_K.
	\end{equation*}
\end{lemma}
\begin{proof}
	Since $\glod$ is piecewise constant and by the definition of $\philod$, we obtain that
	\begin{equation*}
		(\glod, \id_K)_{L^2(K)} = (\glod, b_K)_{L^2(K)}=\ba(\philod, b_K) = \ba((1-\corrc)b_T, b_K).
	\end{equation*}
	Hence, the Fr\'echet derivative of $(\glod, \id_K )_{L^2(K)}$ equals
	\begin{equation*}
		\frac{\partial}{\partial \Ab} (\glod, \id_K )_{L^2(K)} (\delta \Ab) = \frac{\partial}{\partial \Ab} \ba(b_T, b_K) (\delta \Ab) - \frac{\partial}{\partial \Ab} \ba(\corrc b_T, b_K) (\delta \Ab).
	\end{equation*}
	The first term is easily calculated, yielding
	\begin{equation*}
		\frac{\partial \ba(b_T, b_K)}{\partial \Ab} (\delta \Ab) = \int_{D_T} \delta \Ab \nabla b_T \cdot \nabla b_K \dx.
	\end{equation*}
	For the second term, we obtain with the product rule that
	\begin{equation*}
		\frac{\partial\ba(\corrc b_T, b_K)}{\partial \Ab}(\delta \Ab) = \int_{D_T} \delta \Ab \nabla \corrc b_T \cdot \nabla b_K \dx + \int_{D_T}\Ab \nabla \frac{\partial \corrc b_T }{\partial \Ab}(\delta \Ab)\cdot \nabla b_K \dx.
	\end{equation*}
	Using \cref{eq:corrc}, the fact that $\frac{\partial \corrc b_T }{\partial \Ab}(\delta \Ab) \in \mathcal{W}_{T,\ell}$ and the symmetry of $\Ab$, yields that
	\begin{equation}
		\int_{D_T}\Ab \nabla \frac{\partial \corrc b_T }{\partial \Ab}(\delta \Ab)\cdot \nabla b_K \dx = \int_{D_T}\Ab \nabla \frac{\partial \corrc b_T }{\partial \Ab}(\delta \Ab)\cdot \nabla \corrc b_K \dx.
		\label{eq:defcorrc}
	\end{equation} 
	Furthermore, by differentiating \cref{eq:corrc}, we get for any $w \in \mathcal{W}_{T,\ell} $ that
	\begin{align}
		\begin{split}
			\int_{D_T} \delta \Ab \nabla b_T \cdot \nabla w \dx = \int_{D_T} \delta \Ab \nabla \corrc b_T \cdot \nabla w\dx + \int_{D_T }\Ab \nabla \frac{\partial \corrc  b_T }{\partial \Ab}(\delta \Ab)\cdot \nabla w \dx.
		\end{split}
		\label{eq:corrcderivate}
	\end{align}
	Using \eqref{eq:defcorrc} and \cref{eq:corrcderivate} for $w = \corrc b_K $, we obtain for the Fr\'echet derivative that
	\begin{align*}
		\frac{\partial \ba(\corrc b_T, b_K)}{\partial {\Ab}}(\delta {\Ab}) &= \int_{D_T} \delta \Ab \nabla \corrc b_T \cdot \nabla b_K  \dx +  \int_{D_T}\delta \Ab \nabla b_T \cdot \nabla \corrc b_K \dx \\&\quad- \int_{D_T }\delta \Ab \nabla \corrc  b_T \cdot \nabla \corrc b_K \dx.
	\end{align*}
	The $L^2$-representation of the Fr\'echet derivative of $\ba(\corrc b_T, b_K)$ is therefore given by
	\begin{align*}
		\frac{\partial}{\partial \Ab} \ba(\corrc b_T, b_K)  = \nabla \corrc b_T \otimes \nabla b_K + \nabla b_T \otimes \nabla \corrc b_K - \nabla \corrc b_T \otimes \nabla \corrc b_K.
	\end{align*}
	The combination of the above results yields the assertion.
\end{proof}

The following result is needed to estimate the terms appearing after applying the spectral gap inequality in the proof of \cref{lemma:ub}.

\begin{lemma}[$L^4$-regularity estimate for LOD correction operators]\label{l:lemmaregc}
	Let $\Ab$ be a random coefficient field subject to \cref{assumption_coefficient}. Then, the correction of the bubble functions $\corrc b_T$ satisfies the following regularity estimate
	\begin{equation*}
		\int_{D_T} \mathbb{E}\Bigg[\bigg(\fint_{B_\varepsilon(x)}  |\nabla \corrc b_T |^2 \dtx \bigg)^{2} \Bigg] \dx \lesssim \left(\frac{\ell}{H}\right)^{4-d}.
	\end{equation*}
\end{lemma}
\begin{proof}
	First, let $\omega \in \Omega$ be arbitrary but fixed. In order to apply \cref{lem:RegularityAnnealed}, we need to establish the appropriate right-hand side, which results in the equation for $\corrc b_T$ taking the form as in~\cref{e:rhs}. Naturally, $\corrc b_T$ solves, together with the Lagrange multiplier~$\ptl$ the following saddle-point problem 
		\begin{align}\label{eq:saddle-point}
				\begin{pmatrix}
						\Apatch & \BB^T\\\BB & 0
					\end{pmatrix}
				\begin{pmatrix}
						\corrc b_T\\\ptl
					\end{pmatrix} = 
				\begin{pmatrix}
						\Apatch \, b_T\\0
					\end{pmatrix}
			\end{align}
		with the patch-local operators $\Apatch\colon H^1_0(D_T )\rightarrow H^{-1}(D_T), u\mapsto - \nabla \cdot (\Ab \nabla u)$, $\BB\colon H^1_0(D_T )\rightarrow \mathbb P^0(\TT_{H, D_T}),  $ $v\mapsto \Pi_H\vert_{D_T} v$, and its transpose defined by $\BB^T\colon \mathbb P^0(\TT_{H, D_T})\rightarrow H^{-1}(D_T), p \mapsto \{v\in H_0^1(D_T)\mapsto \int_{D_T}p\,v\dx\}$. 
	
	It is a direct consequence that $\corrc b_T$ solves
	\begin{equation*}
		\nabla \cdot (\Ab \nabla \corrc b_T) = \nabla \cdot (\Ab \nabla b_T ) + \BB^T\ptl,
	\end{equation*}
which, for some $\qtl \in L^2(D_T)$, can be rewritten as 
	\begin{equation*}
		\nabla \cdot (\Ab \nabla \corrc b_T) = \nabla \cdot (\Ab \nabla b_T + \qtl).
	\end{equation*}
To see this, we set $\qtl \coloneqq \nabla v $, where $v $ solves $\Delta v = \ptl$ with homogeneous Dirichlet boundary conditions in a ball of radius $C\ell H$, where the constant $C>0$ is chosen such that the ball contains $D_T$.
	
	Furthermore, the local LOD source terms satisfy $\glod = \ptl$; see \cite{HaPe21b}. Hence, using \cref{eq:norm_glod} and following the proofs of \cite[Lem.~4.9]{FisGP19ppt} and \cref{lem:lemmaregphiloc} yields that
		\begin{equation*}
			\label{eq:estgradq}
			\int_{D_T}|\nabla \qtl|^2 \dx =  \int_{D_T}|D^2 v|^2 \dx  \lesssim \int_{D_T}|\ptl|^2 \dx \lesssim H^{d-4},
		\end{equation*}
	as well as $ \normL2{\qtl}{D_T}\lesssim \ell H \normL2{\ptl}{D_T}$, leading to $\int_{D_T}|\qtl|^q\dx \lesssim \ell^{d+q(2-d)/2}H^{d-q}$. Moreover, using $\|\nabla b_T\|_{L^\infty} \approx H^{-1}$ we obtain that $\int_{T}|\Ab\nabla b_T |^q\dx \lesssim H^{d-q} $. Therefore, applying \cref{lem:RegularityAnnealed} for $Q = D_T$, $p=4$ and $\bb = \Ab \nabla b_T - \boldsymbol{q}_{T,\ell}$ yields that
	\begin{align*}
			&\int_{D_T} \mathbb{E}\Bigg[\bigg(\fint_{B_\varepsilon(x)}  |\nabla \corrc b_T |^2 \dtx \bigg)^2 \Bigg] \dx\\
			&\quad\lesssim |D_T|^{(q-4)/q} \bigg(\mathbb E \Bigg[\int_{D_T } |\Ab \nabla b_T |^{q}\dx + \int_{D_T}|\qtl|^{q}\dx\Bigg]\bigg)^{4/q}\\
			&\quad \lesssim \left(\frac{\ell}{H}\right)^{4-d},
		\end{align*}
	which is the assertion.
\end{proof}

\bibliographystyle{alpha}
\bibliography{bib}
\end{document}